\documentclass{amsart}

\usepackage{amsmath,amssymb}

\usepackage{amsfonts}
\usepackage{graphicx}
\usepackage{color}

\def\N{\ifmmode{\mathbb N}\else{$\mathbb N$}\fi} 
\def\R{\ifmmode{\mathbb R}\else{$\mathbb R$}\fi} 
\def\Q{\ifmmode{\mathbb Q}\else{$\mathbb Q$}\fi} 
\def\C{\ifmmode{\mathbb C}\else{$\mathbb C$}\fi} 
\def\Z{\ifmmode{\mathbb Z}\else{$\mathbb Z$}\fi} 

\newtheorem{thm}{Theorem}[section]
\newtheorem{cor}[thm]{Corollary}
\newtheorem{lem}[thm]{Lemma}
\newtheorem{prop}[thm]{Proposition}

\newtheorem{defn}{Definition}[section]

\theoremstyle{remark}
\newtheorem{rem}{\rm\bfseries{Remark}}[section]

\newtheorem{claim}{Claim $(k,i)$}

\begin{document}

\title[Canonical Models and Morse complexes]{Canonical models of filtered $A_{\infty}$-algebras 
and Morse complexes} 

\author[Fukaya, Oh, Ohta and Ono]{Kenji Fukaya \and Yong-Geun Oh \and Hiroshi Ohta \and Kaoru Ono}

\address{K. Fukaya: Department of Mathematics, Kyoto University, 
Kyoto, 606-8502, Japan}
\address{Y.-G. Oh: 
Department of Mathematics, University of Wisconsin, Madison, 
WI 53706, U.S.A. and Korea Institute for Advanced Study, Seoul, Korea}
\address{H. Ohta: Graduate School of Mathematics, Nagoya University, 
Nagoya, 464-8602, Japan}
\address{K. Ono: Department of Mathematics, Hokkaido University, 
Sapporo, 060-0810, Japan}

\thanks{K. F. is partially supported by 
JSPS Grant-in-Aid for Scientific Research No. 18104001, 
Y. O. by US NSF grant  No. 0503954, H. O. by JSPS Grant-in-Aid for 
Scientific Research No. 19340017, and 
K. O. by JSPS Grant-in-Aid for Scientific Research 18340014.}

\maketitle

\centerline{\it In honor of Yasha Eliashberg's sixtieth birthday}

\medskip
\medskip

\section{Introduction.} 

In the book \cite{FOOO}, the authors studied the moduli spaces 
of bordered stable 
maps of genus 0 with Lagrangian boundary condition 
in a systematic way and constructed 
the filtered $A_{\infty}$-algebra associated to Lagrangian 
submanifolds.  
Since our construction depends on various auxiliary choices, 
we considered the {\it canonical model} of filtered 
$A_{\infty}$-algebras, which is unique up to filtered 
$A_{\infty}$-isomorphisms.  
The aim of this note is to explain the construction of the canonical 
model and to apply such an argument to obtain the filtered 
$A_{\infty}$-structure to the Morse complex on the Lagrangian 
submanifold.  
The resulting filtered $A_{\infty}$-operations are described 
by the moduli 
spaces of certain configurations consisting of pseudo-holomorphic 
curves and gradient flow lines.   
Note that the first named author \cite{Fuk94} studied 
the quantization of Morse homotopy based on 
the moduli spaces of certain configurations consisting 
of pseudo-holomorphic discs and gradient flow lines of {\it multiple} 
Morse functions, 
see \cite{Oh2} for monotone case and Theorem A4.28 in \S A 4 in 
\cite{2000}.  
Such configurations are also studied in monotone case by 
Buhovsky \cite{Buk} and Biran and Cornea \cite{BC}.  
We follow Chapter 5 in \cite{FOOO} to explain the algebraic aspect 
of canonical models and use the geometric construction in Chapter 7 
in \cite{FOOO}.  

We briefly review the background of our study.  
Floer \cite{Fl} invented a new theory, which is now called 
Floer (co)homology for Lagrangian intersections.  
Very roughly speaking, it is an analog of Morse theory for 
the action functional on the space of paths with end points 
on Lagrangian submanifolds.  
For a transversal pair of Lagrangian submanifolds $L_0, L_1$, 
the cochain complex is generated by the intersection points 
of $L_0$ and $L_1$.  The coboundary operator is defined by 
counting connecting orbits joining the intersection points.  
The theory was extended by the second named author \cite{Oh} 
to the class of 
monotone Lagrangian submanifolds with the minimal Maslov number 
being at least 3.  
In general, however, there arise obstructions to constructing the Floer 
cochain complex caused by the bubbling-off of pseudo-holomorphic 
discs in the moduli space of connecting orbits.  
We started a systematic study of the moduli spaces of 
pseudo-holomorphic discs with Lagrangian boundary condition 
and formulated the obstructions in terms of the Mauer-Cartan 
equation on the filtered $A_{\infty}$-algebra associated to 
the Lagrangian submanifold \cite{FOOO}. 
In order to give consistent orientations on the moduli spaces, 
we introduced the notion of relative spin structure and 
cosidered relative spin Lagrangian submanifolds.  
For a relative spin pair $(L_0,L_1)$ of Lagrangian submanifolds, 
we constructed a filtered $A_{\infty}$-bimodule over 
the $A_{\infty}$-algbras associated to $L_0$ and $L_1$.  
If each $L_i$, $i=0,1$, admits a solution $b_i$ of 
the Maurer-Cartan equation, we can rectify the Floer operator 
to obtain a coboundary operator $\delta^{b_0,b_1}$.  
Hence the Floer complex 
$(CF^{\bullet}(L_0,L_1), \delta^{b_0,b_1})$ is obtained.  
We also considered the case that the Lagrangian submanifolds 
admit solutions of the Maurer-Cartan equation modulo multiples of 
the fundamental class $[L_i]$ ({\it weak solution}).  
For a weak solution $b_i$, we assign the potential 
${\mathfrak{PO}}(b_i)$. 
If ${\mathfrak{PO}}(b_0)={\mathfrak{PO}}(b_1)$, we can construct 
the Floer complex $(CF(L_0,L_1),\delta^{b_0,b_1})$ deformed by 
$b_0, b_1$.  
This extension with the weak bounding cochains plays crutial role 
in our study of Floer theory on compact toric manifolds \cite{FOOO2}.   

We firstly constucted the filtered $A_{\infty}$-algebra 
mentioned above on suitable subcomplex of 
the singular cochain complex of $L_i$ using systematic 
multi-valued perturbation of Kuranshi maps describing 
the moduli spaces.  
We briefly review these constructions in subsequent section.  
Thus the resulting filtered $A_{\infty}$-algebra depends 
on various choices, i.e., the choice of the subcomplex, 
the choice of systematic multi-valued perturbation, etc.  
In order to make the construction canonical, we introduced 
the notion of the {\it canonical model}.  
Since the structure constants of the filtered $A_{\infty}$-algebra 
depends on these choices, it is appropriate to work with 
the canonical model when we make practical computation of 
the structure constants.  
When we consider ${\mathfrak{PO}}(b)$ as a function 
on the set of weak solutions of the Maurer-Cartan equation, 
we call it the {\it potential function}.  
The canonical model provides an appropriate domain of 
the definition of the potential function.  
The canonical models also play a role in the convergence of a certain 
spectral sequence, see Chapter 6 in \cite{FOOO}.  
(In \cite{2000}, we used another kind of 
finitely generated complex to ensure the weak finiteness property 
of the filtered $A_{\infty}$ algebras.)    
It may be also worth mentioning that we rely on canonical models in 
some places in \cite{FOOO}, since the degree of the ordinary 
cohomology is bounded, though the degree of the singular complex is 
not bounded above.  

We also developed 
an algebraic theory for filtered $A_{\infty}$-algebras, bimodules, 
in particular, 
the homotopy theory of the filtered $A_{\infty}$-algebras, 
bimodules and proved that the homotopy type 
of the resulting algebraic object does not depend on such choices.  
We can also reduce the filtered $A_{\infty}$-structure to 
appropriate free subcomplexes of the original complex.  
In particular, if we work over the ground coefficient field, 
we obtain the filtered $A_{\infty}$-structure on the {\it classical} 
(co)homology of the complex.  
In this note, we review the construction of the canonical models 
of filtered $A_{\infty}$-algebras and filtered $A_{\infty}$-bimodules 
and explain its implication in a geometric setting.  

In section 5, we induce the filtered $A_{\infty}$-algebra structure 
on Morse complex based on the argument in the construction of 
canonical models.  
We choose a Morse function $f$ on $L$, which is adapted to 
a triangulation of $L$ (see section 5).   
\smallskip

\noindent
{\bf Theorem \ref{Morse}}
{\it 
Let $L$ be a relatively spin Lagrangian submanifold in a closed 
symplectic manifold $(M,\omega)$ and $f$ a Morse function on $L$ 
as above.  
Then Morse complex $CM^*(f)\otimes \Lambda_{nov}$ 
carries a structure of a filtered $A_{\infty}$-algebra, 
which is homotopy equivalent to the filtered $A_{\infty}$-algebra 
associated to $L$ constructed in \cite{FOOO}.  
}

\smallskip

This note is based on the lecture at Yashafest and we thank  
the organizers for the invitation.  

\section{Flitered $A_{\infty}$-algebras, bimodules} 
In this subsection, we recall the definition of (filtered) 
$A_{\infty}$-algebras, bimodules, homomorphisms 
and prepare necessary notations.  
Then we explain the notion of homotopy between (filtered) 
$A_{\infty}$-homomorphisms.  
In fact, the notion of homotopy between 
$A_{\infty}$-homomorphisms can be found in the literature, 
e.g., \cite{Smi}.  
For differential graded algebras, homotopy theory was studied 
in rational homotopy theory, see \cite{Su}, \cite{G-M}.  
In order to make clear the relation among such notions, 
we introduced the notion of the model of 
$[0,1] \times \overline{C}$ ($[0,1] \times C$) for 
a (filtered) $A_{\infty}$-algebra $C$ ($\overline{C}$) 
and defined the notion of homotopy using such a model.  

\smallskip
\noindent
{\bf 2.1) Unfiltered $A_{\infty}$-algebras, homomorphisms, bimodules, 
bimodule homomorphisms.}

Let $R$ be a commutative ring, e.g., $\Z$, $\Q$.  
Let $\overline{C}^{\bullet}$ be a cochain complex over $R$.  
We assume that $\overline{C}^k = 0$ for $k < 0$.  
Denote by $\overline{\mathfrak m}_1$ its differential.  
Set $\overline{C}[1]^k=\overline{C}^{k+1}$ and denote the shifted degree by 
$\deg' x = \deg x -1$, where $\deg$ is the original degree 
of $\overline{C}^{\bullet}$.  
In this section, we use only shifted degrees.  
Consider a series of $k$-ary operations, $k=1,2, \dots,$ 

\[
\overline{\mathfrak m}_k : 
(\overline{C}[1]^{\bullet})^{\otimes k} \to \overline{C}[1]^{\bullet}
\]
of degree $1$ with respect to the shifted degrees.  

Before giving the definition of $A_{\infty}$-algebras, we explain 
the case of differential graded algebras.  
Let $(\overline{C}^{\bullet}, d, \cdot)$ be a differential 
graded algebra.  
Define $\overline{\mathfrak m}_1(x)=(-1)^{\deg x} da$ and 
$\overline{\mathfrak m}_2(x \otimes y) = 
(-1)^{\deg x \cdot (\deg y +1)} x \cdot y$.  
Then we find that 

\begin{align}
&\overline{\mathfrak m}_1 \circ \overline{\mathfrak m}_1 (x)= 0,
\nonumber \\
&\overline{\mathfrak m}_1 \circ \overline{\mathfrak m}_2(x \otimes y) 
+ \overline{\mathfrak m}_2 (\overline{\mathfrak m}_1(x_1) \otimes x_2) 
+ (-1)^{\deg' x_1} \overline{\mathfrak m}_2 (x_1 \otimes 
\overline{\mathfrak m}_1(x_2))=0, 
\nonumber \\
&\overline{\mathfrak m}_2(\overline{\mathfrak m}_2(x_1 \otimes x_2) \otimes x_3) + (-1)^{\deg' x_1} \overline{\mathfrak m}_2(x_1 \otimes 
\overline{\mathfrak m}_2(x_2 \otimes x_3))=0, \nonumber
\end{align}
which follow from the facts that $d$ is a differential, 
the multiplication and the differential $d$ satisfies Leibniz' rule 
and the multiplication is associative.  

There are some geometric situations where multiplicative structures 
are defined but not exactly associative.  
A typical example is the composition in based loop spaces.  
In fact, Stasheff \cite{St} introduced the notion of 
$A_{\infty}$-structure on topological spaces in order to characterize 
the homotopy types of based loop spaces.  
He also defined the $A_{\infty}$-structure in algebraic setting.  
For instance, a multiplicative structure is said to be associative 
up to homotopy, if there exists 
$\overline{\mathfrak m}_3: (\overline{C}[1]^{\bullet})^{\otimes 3} \to 
\overline{C}[1]^{\bullet}$ such that 

\begin{align}
& \overline{\mathfrak m}_2(\overline{\mathfrak m}_2(x_1 \otimes x_2) \otimes x_3) + (-1)^{\deg' x_1} \overline{\mathfrak m}_2(x_1 \otimes 
\overline{\mathfrak m}_2(x_2 \otimes x_3)) \nonumber \\ 
+ & \overline{\mathfrak m}_1 \circ \overline{\mathfrak m}_3
(x_1 \otimes x_2 \otimes x_3) + 
\overline{\mathfrak m}_3(\overline{\mathfrak m}_1(x_1) 
\otimes x_2 \otimes x_3) \nonumber \\
+ & (-1)^{\deg' x_1} \overline{\mathfrak m}_3
(x_1 \otimes \overline{\mathfrak m}_1(x_2) \otimes x_3)  
+ (-1)^{\deg' x_1 + \deg' x_2} \overline{\mathfrak m}_3 
(x_1 \otimes x_2 
\otimes \overline{\mathfrak m}_1(x_3)) \nonumber \\
= & 0.  \nonumber
\end{align}

Note that it coincides with the relation corresponding to the 
associativity, if $\overline{\mathfrak m}_3=0$.  
We can continue higher homotopies in a similar way:  
$$
\sum_{k_1+k_2=k+1} \sum_i (-1)^{\sum_{j=1}^{i-1} \deg' x_j} 
\overline{\mathfrak m}_{k_1}
(x_1, \dots, \overline{\mathfrak m}_{k_2}
(x_i, \dots, x_{i+k_2-1}), \dots, x_k) = 0.
$$  
Here $k_1, k_2$ are positive integers.  
For a concise description of relations among higher homotopies, 
we introduce the bar complex of $\overline{C}^{\bullet}$, which is 
defined by 
\[
B(\overline{C}[1]^{\bullet})=\bigoplus_{k=0}^{\infty} B_k(\overline{C}[1]^{\bullet}), \ \ 
B_k(\overline{C}[1]^{\bullet})=\bigoplus_{m_1,\dots, m_k} 
\overline{C}[1]^{m_1} \otimes \cdots \otimes \overline{C}[1]^{m_k}, 
\]
which we consider as a tensor coalgebra.  
The comultiplication is given by 
$$
\Delta (x_1 \otimes \dots \otimes x_k) 
= \sum_{i=0}^k (x_1 \otimes \dots \otimes x_i) \otimes (x_{i+1} \otimes 
\dots \otimes x_k),
$$ 
where $x_1 \otimes \dots \otimes x_i, x_{i+1} \otimes \dots \otimes 
x_k \in B(\overline{C}[1]^{\bullet})$ and 
the former with $i=0$ and the latter with $i=k$ are understood as 
$1 \in B_0(\overline{C}[1]^{\bullet})$.  
Extend $\overline{\mathfrak m}_k$ to the graded coderivation 
$\widehat{\overline{\mathfrak m}_k}$ on $B(\overline{C}[1]^{\bullet})$. 
Namely, 

\begin{align}
& \widehat{\overline{\mathfrak m}}_k(x_1 \otimes \dots \otimes x_N) 
\nonumber \\
= & \sum_{i=1}^{N-k+1} (-1)^{\sum_{j=1}^{i-1} \deg' x_j} 
x_1 \otimes \dots \otimes x_{i-1} \otimes 
\overline{\mathfrak m}_k (x_i \otimes \dots \otimes x_{i+k-1}) \otimes 
x_{i+k} \otimes \dots \nonumber \\
& \hspace{0.5in} \cdots \otimes x_N. \nonumber
\end{align}

We call $(\overline{C}^{\bullet}, \{\overline{\mathfrak m}_k\})$ 
an $A_{\infty}$-algebra, if 
\[
\widehat{\overline{d}}=\sum_k \widehat{\overline{\mathfrak m}}_k:
B(\overline{C}[1]^{\bullet}) \to B(\overline{C}[1]^{\bullet})
\] 
satisfies $\widehat{\overline{d}} \circ \widehat{\overline{d}} = 0$.  
In the case that $\overline{\mathfrak m}_k=0$ for $k > 2$, 
this condition 
is equivalent to the notion of differential graded algebras.  

For a collection $\{\overline{\mathfrak f}_k:B_k(\overline{C}[1]^{\bullet}) \to 
\overline{C}'[1]^{\bullet}\}_{k=1}^{\infty}$ of degree $0$, 
we extend it to a homomorphism as tensor coalgebras
\[
\widehat{\overline{f}}(x_1 \otimes \cdots \otimes x_k)  
= 
\sum_{k_1 + \cdots + k_n=k} 
\overline{\mathfrak f}_{k_1}(x_1 \otimes \cdots \otimes x_{k_1}) 
\otimes \cdots \otimes 
\overline{\mathfrak f}_{k_n}(x_{k+1-k_n} \otimes \cdots \otimes x_k).
\]
We call $\{\overline{\mathfrak f}_k\}$ an $A_{\infty}$-homomorphism, 
if $\widehat{\overline{f}}$ satisfies 
$\widehat{\overline{d}}_{\overline{C}'} \circ \widehat{\overline f} = 
\widehat{\overline f} \circ \widehat{\overline{d}}_{\overline{C}}.$  

In terms of the components $\overline{\mathfrak m}_k$'s and 
$\overline{\mathfrak f}_k$'s, 
this is equivalent to 

\begin{align} 
& \sum_{i_1+ \dots + i_k=n} \overline{\mathfrak m}_k(
\overline{\mathfrak f}_{i_1}(x_1 \otimes \dots \otimes x_{i_1}) 
\otimes \dots 
\otimes \overline{\mathfrak f}_{i_k}(x_{i_1+ \dots + i_{k-1} +1} 
\otimes 
\dots \otimes x_{n})) \nonumber \\
= & 
\sum_{j_1+ \dots + j_{\ell}=n} \sum_{p=1}^{\ell} 
(-1)^{\sum_{i=1}^{j_1+ \dots j_{p-1}} \deg' x_i}
\overline{\mathfrak f}_{j_1}(x_1 \otimes \dots \otimes x_{j_1}) \otimes 
\dots \otimes \nonumber \\
& \overline{\mathfrak m}_{j_p}(x_{j_1+ \dots + j_{p-1}+1} 
\otimes \dots \otimes 
x_{j_1+ \dots + j_p}) \otimes \dots \otimes 
\overline{\mathfrak f}_{j_{\ell}}(x_{j_1+ \dots j_{\ell -1} +1} 
\otimes \dots \otimes x_n) \nonumber
\end{align}

Let $(\overline{C}^{\bullet}_i, \{\overline{\mathfrak m}_k^{(i)}\})$, 
$i=0,1$, 
be $A_{\infty}$-algebras, $\overline{D}^{\bullet}$ a graded module and 
$\overline{\mathfrak n}_{k_1,k_0}:B_{k_1}(\overline{C}_1[1]^{\bullet}) \otimes 
\overline{D}[1]^{\bullet} \otimes B_{k_0}(\overline{C}_0[1]^{\bullet}) \to 
\overline{D}[1]^{\bullet}$ homomorphisms of degree $1$.  
We call $(\overline{D}^{\bullet},\{\overline{\mathfrak n}_{k_1,k_0}\})$ 
an $A_{\infty}$-bimodule, if 
\[
\widehat{\overline{d}}_{\overline{\mathfrak n}} \circ 
\widehat{\overline{d}}_{\overline{\mathfrak n}} = 0, 
\]
where $\widehat{\overline{d}}_{\overline{\mathfrak n}}$ is 
defined on 
$
B(\overline{C}_1[1]^{\bullet}) \otimes \overline{D}[1]^{\bullet} \otimes 
B(\overline{C}_0[1]^{\bullet})
$ as follows:  

\begin{align}
& \widehat{\overline{d}}_{\overline{\mathfrak n}}
(x_{1,1} \otimes \dots \otimes x_{1,k_0} \otimes y \otimes 
x_{0,1} \otimes \dots \otimes x_{0,k_0}) \nonumber \\
= & \widehat{\overline{d}}^{(1)}(x_{1,1} \otimes \dots 
\otimes x_{1,k_1}) 
\otimes y \otimes x_{0,1} \otimes \dots \otimes x_{0,k_0} \nonumber \\
& + \sum_{k'_1 \leq k_1, k'_0 \leq k_0} 
(-1)^{\sum_{i=1}^{k_1-k'_1} \deg' x_i}
x_{1,1} \otimes \dots x_{k_1-k'_1} \nonumber \\
& \otimes {\overline{\mathfrak n}}_{k'_1,k'_0}(x_{1,k_1-k'_1+1} 
\otimes \dots x_{k_1} \otimes y \otimes x_{0,1} 
\otimes \dots 
\otimes x_{0,k'_0}) \otimes x_{0,k_0+1} \otimes \dots \otimes x_{0,k_0} 
\nonumber \\
& + (-1)^{\sum_{i=1}^{k_1} \deg' x_i + \deg' y}
x_{1,1} \otimes \dots x_{1,k_1} \otimes y \otimes \widehat{\overline{d}}^{(0)}(x_{0,1} \otimes \dots \otimes x_{0,k_0}). \nonumber
\end{align}

The condition for $(\overline{D},\{\overline{\mathfrak n}_{k_1,k_0}\})$ 
to be an $A_{\infty}$-bimodules over $\overline{C}_i$, $i=0,1$ is 
equivalent to the identity 
 
\begin{align}
& \overline{\mathfrak n}_{*,*}(\widehat{\overline{d}}^{(1)}(
x_{1,1} \otimes \dots \otimes x_{1,k_1}) \otimes y \otimes 
x_{0,1} \otimes \dots \otimes x_{0,k_0}) \nonumber \\
+ & \sum_{k'_1 \leq k_1, k'_0 \leq k_0} (-1)^{\sum_{i=1}^{k_1-k'_1} 
\deg' x_i} 
\overline{\mathfrak n}_{k_1-k'_1,k_0-k'_0}
(x_{1,1} \otimes \dots \otimes 
x_{k_1-k'_1} \otimes \nonumber \\
& \overline{\mathfrak n}_{k'_1,k'_0}
(x_{1,k_1-k'_1+1} \otimes \dots \otimes x_{1,k_1} \otimes y 
\otimes x_{0,1} \otimes x_{0,k'_0}) \otimes x_{0,k'_0+1} \otimes 
\dots \otimes x_{0,k_0}) \nonumber \\
+ & (-1)^{\sum_{i=1}^{k_1} \deg' x_i + \deg' y} 
\overline{\mathfrak n}_{*,*} (x_{1,1} \otimes \dots \otimes 
x_{1,k_1} \otimes y \otimes \widehat{\overline{d}}^{(0)} 
(x_{0,1} \otimes \dots \otimes x_{0,k_0})) \nonumber \\
= & 0. \nonumber
\end{align}

Here $\overline{\mathfrak n}_{*,*}:B(\overline{C}_1[1]^{\bullet}) 
\otimes \overline{D}[1] \otimes B(\overline{C}_0[1]^{\bullet}) 
\to \overline{D}[1]$ is defined to be 
$\overline{\mathfrak n}_{k_1,k_0}$ on 
$B_{k_1}(\overline{C}_1[1]^{\bullet}) 
\otimes \overline{D}[1] \otimes B_{k_0}(\overline{C}_0[1]^{\bullet})$.  

Let $\{\overline{\mathfrak f}_k^{(i)}:\overline{C}_i[1]^{\bullet} \to 
\overline{C}_i'[1]^{\bullet}\}$, $i=0,1$, be $A_{\infty}$-homomorphisms 
and $\overline{D}^{\bullet}$, resp. $\overline{D'}^{\bullet}$ 
an $A_{\infty}$-bimodules over 
$\overline{C}_i$, resp. $\overline{C}_i'$, $i=0,1$.  
For a collection $\{\overline{\phi}_{k_1,k_0}\}:
B_{k_1}(\overline{C}_1[1]^{\bullet}) \otimes \overline{D}[1]^{\bullet} 
\otimes 
B_{k_0}(\overline{C}_0[1]^{\bullet}) 
\to \overline{D}'[1]^{\bullet}$ of degree $0$, we define 
\[
\widehat{\overline{\phi}}:B(\overline{C}_1[1]^{\bullet}) \otimes \overline{D}[1]^{\bullet} 
\otimes B(\overline{C}_0[1]^{\bullet}) \to 
B(\overline{C}_1'[1]^{\bullet}) \otimes \overline{D}'[1]^{\bullet} \otimes 
B(\overline{C}_0'[1]^{\bullet})
\]
as the homomorphism determined by $\{\overline{\mathfrak f}_k^{(i)}\}$, 
$i=0,1$ and $\{\overline{\phi}_{k_1,k_0}\}$.  
We call $\{\overline{\phi}_{k_1,k_0}\}$ a homomorphism of 
$A_{\infty}$-bimodules, if 
\[
\widehat{\overline{d}}_{\overline{\mathfrak n}'} \circ 
\widehat{\overline{\phi}} = \widehat{\overline{\phi}} \circ 
\widehat{\overline{d}}_{\overline{\mathfrak n}}.
\]
 
\smallskip
\noindent 
{\bf 2.2) The universal Novikov ring and the energy filtration.}

To explain the notion of filtered $A_{\infty}$-algebras,   
we introduce the universal Novikov ring.  
Let $e$ and $T$ be formal variables of degree $2$ and $0$, 
respectively.  
Set  
\begin{align}
 \Lambda_{nov} & = \{ \sum_i a_i e^{\mu_i} T^{\lambda_i} \ 
\vert \ a_i \in R, 
\ \mu_i \in \Z, \ \lambda_i \in \R, \ 
\lambda_i \to +\infty (i \to +\infty)\} \nonumber \\
 \Lambda_{0,nov} & = 
\{ \sum_i a_i e^{\mu_i} T^{\lambda_i} \in \Lambda_{nov} 
\ \vert \ \lambda_i \geq 0 \}. \nonumber 
\end{align}

Set 
$$C^{\bullet}=
\{\sum c_i e^{\mu_i} T^{\lambda_i} \vert 
c_i \in \overline{C}^{\bullet}, \mu_i \in \Z, \lambda_i \in \R, 
\lambda_i \to +\infty (i \to +\infty)\},
$$
which is the completion of the graded tensor product 
$\overline{C}^{\bullet} \otimes_R \Lambda_{0,nov}$ with respect to 
the energy filtration given below. 
We define the filtration defined by 
\[
F^{\lambda}C^{\bullet}= \{ \sum_i x_i e^{\mu_i} T^{\lambda_i} \in C^{\bullet}
\ \vert \ x_i \in \overline{C}^{\bullet}, \lambda_i \geq \lambda \} 
\]
on $C^{\bullet}$ and denote by 
$
F^{\lambda}(C[1]^{m_1} \otimes \cdots \otimes 
C[1]^{m_k})
$ 
the submodule of $\C[1]^{m_1} \otimes \cdots \otimes C[1]^{m_k}$ 
spanned by 
\[
F^{\lambda_1}(C[1]^{m_1}) \otimes \cdots \otimes 
F^{\lambda_k}(C[1]^{m_k}), \ \ 
\sum_{i=1}^k \lambda_i = \lambda.
\]  
Define the bar complex of $C[1]^{\bullet}$ by 
the completion with respect to the energy filtration and denote it by  
$B_k(C[1]^{\bullet})$.

\smallskip
\noindent
{\bf 2.3) Filtered case and $G$-gapped conditions.}

Consider the $k$-ary operations, $k=0,1,2, \dots$,  
\[
{\mathfrak m}_k:B_k(C[1]^{\bullet}) \to C[1]^{\bullet}
\]
such that 
\[
{\mathfrak m}_k(F^{\lambda_1}C[1]^{\bullet} \otimes \cdots 
\otimes F^{\lambda_k} C[1]^{\bullet}) \subset F^{\lambda_1 + \cdots 
\lambda_k} C[1]^{\bullet}
\]
and 
\[
{\mathfrak m}_0(1) \in F^{\lambda'}C[1]^{\bullet} 
\text{ for some } \lambda'>0.
\]

We used the induction on the energy level in various arguments 
in \cite{FOOO}, see also section 3 in this note.  
For such purposes, we introduced the 
$G$-gapped condition, which we assume from now on, as follows.  
Note that the $G$-gapped condition follows from Gromov's compactness 
theorem in the case of symplectic Floer theory.  
Let $G \subset \R_{\geq 0} \times 2\Z$ be a monoid such that 
${\rm pr}_1^{-1}([0,c])$ is finite for any $c \geq 0$ and 
\[
{\rm pr}_1^{-1} (0) = \{{\mathbf 0}=(0,0)\}.
\]
Here ${\rm pr}_i$ is the projection to the $i$-th factor, $i=1,2$.  
The filtered $A_{\infty}$-algebra is said to be $G$-gapped, if 
there exist 
\[
{\mathfrak m}_{k,\beta_i}:B_k(\overline{C}[1]^{\bullet}) \to 
\overline{C}[1]^{\bullet}
\]
for $\beta_i=(\lambda_i,\mu_i) \in G \subset \R_{\geq 0} \times 2\Z$ 
such that 
${\mathfrak m}_{0,{\mathbf 0}}=0$ and 
\[
{\mathfrak m}_k = \sum_{i} 
T^{\lambda_i} e^{\mu_i/2} {\mathfrak m}_{k, \beta_i}.  
\] 
Extend ${\mathfrak m}_k$ to the graded coderivation 
$\widehat{\mathfrak m}_k$ on $B(C[1]^{\bullet})$.  
We call $(C^{\bullet},\{{\mathfrak m}_k\})$ a filtered 
$A_{\infty}$-algebra, if 
\[
\widehat{d} = \sum_k \widehat{\mathfrak m}_k:B(C[1]^{\bullet}) \to B(C[1]^{\bullet})
\] 
satisfies $\widehat{d} \circ \widehat{d} = 0$.  
In other words, 
$$
\sum_{k_1+k_2=k+1} \sum_i (-1)^{\sum_{j=1}^{i-1} \deg' x_j} 
{\mathfrak m}_{k_1}
(x_1, \dots, {\mathfrak m}_{k_2}
(x_i, \dots, x_{i+k_2-1}), \dots, x_k) = 0.
$$  
Here $k_1$ is a positive integer and $k_2$ is a non-negative 
integers.  
When $k_2=0$, ${\mathfrak m}_{k_2}(x_i, \dots, x_{i+k_2-1})$ 
is understood as ${\mathfrak m}_0(1)$.  

For a filtered $A_{\infty}$-algebra 
$(C^{\bullet},\{{\mathfrak m}_k\})$,
set $\overline{\mathfrak m}_k= {\mathfrak m}_{k,{\bf 0}}$, 
${\bf 0}=(0,0) \in \R_{\geq 0} \times 2\Z$.  
Then $(\overline{C}^{\bullet},\{\overline{\mathfrak m}_k\})$ is an 
$A_{\infty}$-algebra.  
We call $(C^{\bullet},\{{\mathfrak m}_k\})$ a deformation of 
$(\overline{C}^{\bullet},\{\overline{\mathfrak m}_k\})$.  

Note that ${\mathfrak m}_1 \circ {\mathfrak m}_1$ may not be zero 
and we have 
$$
{\mathfrak m}_1 \circ {\mathfrak m}_1 (x) + {\mathfrak m}_2(
{\mathfrak m}_0(1), x) + (-1)^{\deg' x} {\mathfrak m}_2(x,
{\mathfrak m}_0(1)) = 0.
$$
We set 
$$
e^b=1+b+b\otimes b + b \otimes b \otimes b + \dots, 
$$
for $b \in {\mathcal F}^{\lambda}(C[1]^{0})$ with $\lambda > 0$ 
and consider the {\it Maurer-Cartan} equation: 
$$
\widehat{d}(e^b) = 0,
$$
which is equivalent to 
$$
{\mathfrak m}_0(1)+{\mathfrak m}_1(b)+{\mathfrak m}_2(b,b)+ 
{\mathfrak m}_3(b,b,b) + \dots = 0.
$$
For a given $b$, we define a coalgebra homomorphism
$$
\Phi^b(x_1 \otimes x_2 \otimes \dots \otimes x_k) 
=e^b \otimes x_1 \otimes e^b \otimes x_2 \otimes e^b 
\otimes \dots \otimes e^b \otimes x_k \otimes e^b.
$$
Then define 
$$
{\mathfrak m}_k^b(x_1 \otimes \dots \otimes x_k)
={\mathfrak m}_* \circ \Phi^b(x_1 \otimes \dots \otimes x_k),
$$
where 
${\mathfrak m}_*:B(C[1]^{\bullet}) \to C[1]^{\bullet}$ is defined 
by ${\mathfrak m}_*\vert_{B_k(C[1]^{\bullet})}={\mathfrak m}_k$.  
Then, for a solution $b$ of the Maurer-Cartan equation, we find that 
${\mathfrak m}_0^b(1)=0$, hence ${\mathfrak m}_1^b \circ 
{\mathfrak m}_1^b=0$.  
Namely, the original ${\mathfrak m}_1$ is rectified to 
a coboundary operator ${\mathfrak m}_1^b$ using a solution of 
the Maurer-Cartan equation, which we also call a bounding cochain.

For a collection $\{ {\mathfrak f}_k:B_k(C[1]^{\bullet}) \to C'[1]^{\bullet} \}_{k=0}
^{\infty}$ of degree $0$, we define 
\[
\widehat{\mathfrak f}(x_1 \otimes \cdots \otimes x_k) 
=  
\sum_{k_1 + \cdots k_n=k} 
{\mathfrak f}_{k_1}(x_1 \otimes \cdots \otimes x_{k_1}) 
\otimes \cdots \otimes 
{\mathfrak f}_{k_n}(x_{k+1-k_n} \otimes \cdots \otimes x_k), 
\]
for $k>0$ and 
\[
\widehat{\mathfrak f}(1)=1 + {\mathfrak f}_0(1) + 
{\mathfrak f}_0(1) \otimes {\mathfrak f}_0(1) + \cdots,
\]
where $1 \in \Lambda_{0,nov}=B_0(C[1]^{\bullet})$.  
We assume the $G$-gapped condition, i.e., 
there exist 
\[
{\mathfrak f}_{k,\beta_i}:B_k(\overline{C}[1]^{\bullet}) \to 
\overline{C}'[1]^{\bullet}
\] 
for $\beta_i=(\lambda_i,\mu_i) \in G$ 
with $\lambda_i \to +\infty$ as 
$i \to +\infty$ 
such that 
\[
{\mathfrak f}_k = \sum_i T^{\lambda_i} e^{\mu_i/2} 
{\mathfrak f}_{k,\beta_i}.
\]
In particular, $\widehat{\mathfrak f}$ preserves the energy 
filtration.  Namely, 
\[
\widehat{\mathfrak f}(F^{\lambda}B(C[1]^{\bullet})) \subset F^{\lambda}C'[1]^{\bullet},
\]
where $\{F^{\lambda}B(C[1]^{\bullet})\}$ is the filtration derived from 
the filtration $F^{\lambda}$ on $C[1]^{\bullet}$.  
We call $\{ {\mathfrak f}_k \}$ a $G$-gapped 
filtered $A_{\infty}$-homomorphism, 
if 
$\widehat{d}_{C}' \circ \widehat{\mathfrak f} = \widehat{\mathfrak f} 
\circ \widehat{d}_{C}$.  
When we do not specify the monoid $G$, we call gapped filtered $A_{\infty}$-algebras, gapped filtered $A_{\infty}$-homomorphisms, etc.  

\smallskip
\noindent
{\bf 2.4) Homotopy theory.}

In \cite{FOOO}, we introduced the notion of models of 
$[0,1] \times C^{\bullet}$ (Definition 15.1) and gave two constructions. Using this notion, we developed 
the homotopy theory of filtered
$A_{\infty}$-algebras and filtered-$A_{\infty}$ bimodules.
Our formulation has an advantage to clarify 
equivalence of various definitions of homotopy
of $A_{\infty}$ algebras appearing in the literature
even for the unfiltered cases.

Let $C$ be the completion of $\overline{C} \otimes \Lambda_0$, which is 
a filtered $A_{\infty}$-algebra. 

\begin{defn}
Let ${\mathfrak C}$ be the completion of $\overline{\mathfrak C} \otimes \Lambda_0$, which is 
a filtered $A_{\infty}$-algebra together with 
filtered $A_{\infty}$-homomorphisms.  
\[
{\rm Incl}:C \to {\mathfrak C}, \ 
{\rm Eval}_{s=i}:{\mathfrak C} \to C, i=0,1.
\]  
We call $\mathfrak C$ a model of $[0,1] \times C$, if 
the following conditions are satisfied:  

\begin{itemize}
\item
${\rm Incl}_{k,\beta}$ and ${\rm Eval}_{s=i}, i=0,1$ are 
zero unless $(k,\beta)=(1,\beta_0)$.  
\item
${\rm Incl}_{1,\beta_0}, ({\rm Eval}_{s=0})_{1,\beta_0}$ 
are cochain homotopy equivalences between 
$\overline{C}$ and $\overline{\mathfrak C}$.  
\item
${\rm Eval}_{s=0} \circ {\rm Incl} = {\rm Eval}_{s=1} \circ {\rm Incl} 
= {\rm id.}$  
\item
${\rm Eval}_{s=0} \oplus {\rm Eval}_{s=1}:{\mathfrak C} \to C \oplus C$ 
is surjective.  
\end{itemize}
\end{defn}

We quote here one of constructions of the models of $[0,1] \times C$ 
for reader's convenience.  

Set $$C^{[0,1]}=C \oplus C[-1] \oplus C,$$
and define ${\mathfrak I}_0, {\mathfrak I}_1:C \to C^{[0,1]}$ of 
degree 0 and 
${\mathfrak I}:C \to C^{[0,1]}$ of degree 1 by 
$$
{\mathfrak I}_0(x)=(x,0,0), {\mathfrak I}_1(x)=(0,0,x), 
{\mathfrak I}_1(x)=(0,x,0).
$$
We extend ${\mathfrak I}_0,{\mathfrak I}_1$ to 
$B(C[1]) \to B(C^{[0,1]}[1])$ and denote them by the same symbol.  
Define 

\begin{align}
({\rm Eval}_{s=0})_1(x,y,z)=x, \ ({\rm Eval}_{s=1})_1(x,y,z)=z, 
\nonumber \\
({\rm Incl})_1(x)={\mathfrak I}_0(x)+{\mathfrak I}_1(x)=(x,0,x). 
\nonumber
\end{align}

We define the filtered $A_{\infty}$-structure $\{{\mathfrak M}_k\}$.  

For ${\mathfrak M}_0$, ${\mathfrak M}_1$, we set  
\begin{align}
{\mathfrak M}_0(1) & 
= ({\rm Incl})_1({\mathfrak m}_0(1)), \nonumber \\
{\mathfrak M}_1({\mathfrak I}_0(x)) & 
={\mathfrak I}_0({\mathfrak m}_1(x)) 
+ (-1)^{\deg' x} {\mathfrak I}(x), \nonumber \\
{\mathfrak M}_1({\mathfrak I}_1(x)) & 
={\mathfrak I}_1({\mathfrak m}_1(x)) - (-1)^{\deg' x} {\mathfrak I}(x), 
\nonumber \\
{\mathfrak M}_1({\mathfrak I}(x)) & 
={\mathfrak I}({\mathfrak m}_1(x)).
\nonumber
\end{align}

We define ${\mathfrak M}_k$, $k \geq 2$ as follows.  
For ${\mathbf x} \in B_k(C[1])$, $y \in C$ and 
${\mathbf z} \in B_{\ell}(C[1])$, we set 
$$
{\mathfrak M}_{k+\ell + 1}({\mathfrak I}_0({\mathbf x}), {\mathfrak I}(y), {\mathfrak I}_1({\mathbf z})) = 
(-1)^{\deg' {\mathbf z}} {\mathfrak I}({\mathfrak m}_{k + \ell + 1} 
({\mathbf x}, y, {\mathbf z}),
$$
and 
$$
{\mathfrak M}_k({\mathfrak I}_0({\mathbf x}))=
{\mathfrak I}_0({\mathfrak m}_k({\mathbf x})), 
{\mathfrak M}_{\ell}({\mathfrak I}_1({\mathbf z}))=
{\mathfrak I}_1({\mathfrak m}_{\ell}({\mathbf z})), \ \ 
{\text{ for }} k, \ell \geq 2.
$$  
Here the order of 
${\mathfrak I}_0({\mathbf x}), {\mathfrak I}(y), 
{\mathfrak I}_1({\mathbf z})$ is important.  
We define operators ${\mathfrak M}_k$ on $C^{[0,1]}$ other than 
those defined above to be zero.

Models of $[0,1] \times C$ are not unique, but we proved the following:

\begin{thm}[Theorem 15.34 in \cite{FOOO}]
Let $C_1,C_2$ be gapped filtered $A_{\infty}$-algebras and 
${\mathfrak C}_1, {\mathfrak C}_2$ any models for $[0,1] \times C_1, 
[0,1] \times C_2$, respectively.  
Let ${\mathfrak f}:C_1 \to C_2$ be a gapped filtered 
$A_{\infty}$-homomorphism.  
Then there exists a gapped filtered $A_{\infty}$-homomorphism 
${\mathfrak F}:{\mathfrak C}_1 \to {\mathfrak C}_2$ such that 
$$
{\rm Eval}_{s=s_0} \circ {\mathfrak F} = {\mathfrak f} \circ 
{\rm Eval}_{s=s_0}, \ s_0=0,1
$$
and 
$$
{\rm Incl} \circ {\mathfrak f} = {\mathfrak F} \circ {\rm Incl}.
$$
\end{thm}

We define two filtered $A_{\infty}$-homomorphisms 
${\mathfrak f}_i:C_1 \to C_2, i=0,1$ are homotopic, if 
there is a model ${\mathfrak C}_2$ of $C_2$ and 
a filtered $A_{\infty}$-homomorphism 
${\mathfrak F}:C_1 \to {\mathfrak C}_2$ such that 
${\mathfrak f}_i={\rm Eval}_{s=i} \circ {\mathfrak F}$.  
Although the definition literally depends on the choice of the model 
${\mathfrak C}_2$, 
we can show that the notion of homotopy between ${\mathfrak f}_i$ 
does not depend on the choice of the model and the homotopy 
is, in fact, an equivalence relation, see Chapter 4 in \cite{FOOO}.  

Note that the notion of homotopy between $A_{\infty}$-homomorphisms 
in the unfiltered case appeared in literaturure, e.g., \cite{Smi}.  
By taking a suitable model, we can find that 
our definition above coincides with such a definition.  
It also implies that various definitions which appear in the literature 
are equivalent to one another.  
We think that the notion of models clarifies arguments and is 
also useful when we consider the gauge equivalence between 
solutions of the Maurer-Cartan equation.  
Namely, two solutions $b, b'$ of the Maurer-Cartan equation is 
gauge equivalent, if there is a model ${\mathfrak C}$ of 
$[0,1] \times C$ and a solution $\widetilde{b}$ of the Maurer-Cartan 
eqution on ${\mathfrak C}$ such that 
${\rm Eval}_{s=0}(\widetilde{b})=b$ and 
${\rm Eval}_{s=1}(\widetilde{b})=b'$.  
For details, see section 16 in \cite{FOOO}.  

Among other things, we proved the Whitehead type theorem as follows.  
An $A_{\infty}$-homomorphism $\{\overline{\mathfrak f}_k\}$ from 
$\overline{C}^{\bullet}$ to $\overline{C}^{' \bullet}$ 
is called a weak homotopy equivalence, if 
$\overline{f}_1:\overline{C}^{\bullet} \to \overline{C}^{' \bullet}$ 
is a cochain homotopy equivalence 
between $\overline{\mathfrak m}_1$-complexes. 
A filtered $A_{\infty}$-homomorphism $\{{\mathfrak f}_k\}$ from 
$C[1]^{\bullet}$ to $C'[1]^{\bullet}$ is called 
a weak homotopy equivalence, if 
$\overline{\mathfrak f}_1={\mathfrak f}_{1,{\bf 0}}$ is a cochain 
homotopy equivalence between 
$\overline{\mathfrak m}_1={\mathfrak m}_{1,{\bf 0}}$-complexes.  

\begin{thm}[Theorem 15.45 in \cite{FOOO}]\label{whitehead}
(1) A weak homotopy equivalence of $A_{\infty}$-algebras is 
a homotopy equivalence.  

\noindent
(2) A gapped weak homotopy equivalence between gapped filtered 
$A_{\infty}$-algebras is a homotopy equivalence.  
The homotopy inverse of a strict weak homotopy equivalence can be taken 
to be strict.  
\end{thm}

Note that the above theorem does not hold in the realm of 
differential graded algebras.  The notion of $A_{\infty}$-homomorphism 
is much wider than that of homomorphisms as differential graded 
algebras.  

\smallskip
\noindent
{\bf 2.5) Filtered $A_{\infty}$-bimodules.} 

Let $(C_i^{\bullet},\{{\mathfrak m}_k^{(i)}\})$, $i=0,1$, be filtered 
$A_{\infty}$-algebras and $\overline{D}^{\bullet}$ a graded module.  
Write 
\[
D[1]^{\bullet}=\overline{D}[1]^{\bullet} \otimes \Lambda_{0,nov}
\]
and 
\[
\widetilde{D}[1]^{\bullet}=\overline{D}[1]^{\bullet} \otimes \Lambda_{nov}.
\]
Let 
${\mathfrak n}_{k_1,k_0}:B_{k_1}(C_1[1]^{\bullet}) \otimes D[1]^{\bullet} \otimes 
B_{k_0}(C_0[1]^{\bullet}) \to D[1]^{\bullet}$, $k_1,k_0=0,1,2,\dots$, 
be $\Lambda_{0,nov}$-module homomorphisms of degree $1$.  
We also denote its extension to $\widetilde{D}[1]^{\bullet}$ by the same symbol 
${\mathfrak n}_{k_1,k_0}$.  
We call $(D^{\bullet},\{{\mathfrak n}_{k_1,k_0}\})$ and 
$(\widetilde{D}^{\bullet},\{{\mathfrak n}_{k_1,k_0}\})$ a filtered 
$A_{\infty}$-bimodule over $(C_i^{\bullet},\{{\mathfrak m}_k\})$, if 
\[
\widehat{d}_{\mathfrak n} \circ \widehat{d}_{\mathfrak n} = 0,
\]
where $\widehat{d}_{\mathfrak n}$ is the coderivation 
on $B(C_1[1]^{\bullet}) \otimes D[1]^{\bullet} \otimes B(C_0[1]^{\bullet})$ determined by 
$\{{\mathfrak m}_k^{(i)}\}$, $i=0,1$, and 
$\{{\mathfrak n}_{k_1,k_0}\}$.  
The $G$-gapped condition is defined in a similar way to 
the case of filtered $A_{\infty}$-algebras: 
\[
{\mathfrak n}_{k_1,k_0} = \sum_{\beta \in G} T^{{\rm pr}_1(\beta)} 
e^{{\rm pr}_2(\beta)/2} {\mathfrak n}_{k_1,k_0,\beta}, 
\]
where 
\[
{\mathfrak n}_{k_1,k_0,\beta}:B_{k_1}(\overline{C}_1[1]^{\bullet}) \otimes 
\overline{D}[1]^{\bullet} \otimes 
B_{k_0}(\overline{C}_0[1]^{\bullet}) \to \overline{D}[1]^{\bullet}.  
\]
For $\lambda \in \R$, we set 
\[
{\mathcal F}^{\lambda}\widetilde{D}^{\bullet}=T^{\lambda} \cdot 
D^{\bullet}.
\]  
For $\lambda \geq 0$, 
${\mathcal F}^{\lambda}\widetilde{D}^{\bullet} \subset D^{\bullet}$, 
hence we obtain  
a filtration on $D^{\bullet}$.  
It is clear that 
\[
\widehat{d}_{\mathfrak n}(F^{\lambda_1}B(C_1[1]^{\bullet}) \otimes 
{\mathcal F}^{\lambda_2}\widetilde{D}^{\bullet} 
\otimes F^{\lambda_3}B(C_0[1]^{\bullet})) 
\subset {\mathcal F}^{\lambda_1+\lambda_2+\lambda_3}\widetilde{D}[1]^{\bullet}.  
\]

Note that ${\mathfrak n}_{0,0} \circ {\mathfrak n}_{0,0}$ may not 
be zero and we have 
$$
{\mathfrak n}_{0,0} \circ {\mathfrak n}_{0,0}(y) 
+ {\mathfrak n}_{1,0}({\mathfrak m}_0^{(1)}(1), y) 
+(-1)^{\deg' y} {\mathfrak n}_{0,1}(y,{\mathfrak m}_0^{(0)}(1)) = 0.
$$
For $b_i \in {\mathfrak F}^{\lambda^{(i)}}(C_i[1]^{\bullet})$ 
with $\lambda^{(i)} > 0$, we define 
$$
{\mathfrak n}^{b_0,b_1}_{k_1,k_0}({\mathbf x} \otimes y \otimes 
{\mathbf z})={\mathfrak n}_{*,*}(\Phi^{b_1}({\mathbf x}) 
\otimes y \otimes \Phi^{b_0}({\mathbf z})), 
$$
for ${\mathbf x} \in B_{k_1}(C_1[1]^{\bullet})$ and 
${\mathbf z} \in B_{k_0}(C_0[1]^{\bullet})$.  
In particular, 
$$
{\mathfrak n}_{0,0}^{b_0,b_1}(y)=
{\mathfrak n}_{*,*}(e^{b_1} \otimes y \otimes e^{b_0}).
$$
If $b_0, b_1$ are solutions of the Maurer-Cartan equations 
in the filtered $A_{\infty}$-algebras $C_0, C_1$, respectively, 
we find that 
$${\mathfrak n}_{0,0}^{b_0,b_1} \circ 
{\mathfrak n}_{0,0}^{b_0,b_1} = 0.$$  
Namely, we can rectify the original ${\mathfrak n}_{0,0}$ 
to a coboundary operator ${\mathfrak n}_{0,0}^{b_0,b_1}$ on $D[1]$.  

Let $\{{\mathfrak f}_k^{(i)}\}$, $i=0,1$, be 
filtered $A_{\infty}$-homomorphisms from $C_i^{\bullet}$ to 
$C_i^{' \bullet}$ and 
$(\widetilde{D}^{\bullet},\{{\mathfrak n}_{k_1,k_0}\})$, resp. 
$(\widetilde{D}^{' \bullet},\{{\mathfrak n}_{k_1,k_0}'\})$, 
filtered $A_{\infty}$-bimodules over 
$(C_i^{\bullet},\{{\mathfrak m}_k^{(i)}\})$, resp. 
$(C_i^{' \bullet},\{{\mathfrak m}_k{'(i)}\})$.  
Suppose that there exist a real number $c$ and 
$\Lambda_{nov}$-homomorphisms 
${\phi}_{k_1,k_0}:B_{k_1}(C_1[1]^{\bullet}) \otimes \widetilde{D}[1]^{\bullet} 
\otimes B_{k_0}(C_0[1]^{\bullet}) 
\to \widetilde{D}'[1]^{\bullet}$, $k_1,k_0=0,1,2,\dots$, 
such that 
\[
{\phi}_{k_1,k_0}(F^{\lambda_1}B_{k_1}(C_1[1]^{\bullet}) \otimes 
{\mathcal F}^{\lambda_2} \widetilde{D}[1]^{\bullet} \otimes 
F^{\lambda_3}B_{k_0}(C_0[1]^{\bullet})) \subset 
{\mathcal F}^{\lambda_1+\lambda_2+\lambda_3 -c}\widetilde{D}'[1]^{\bullet}.
\]
We call such $c$ the energy loss of $\{{\phi}_{k_1,k_0}\}$.  
For such a collection $\{{\phi}_{k_1,k_0}\}$, 
we define 
\[
\widehat{\mathfrak \phi}:B(C_1[1]^{\bullet}) \otimes \widetilde{D}[1]^{\bullet} \otimes 
B(C_0[1]^{\bullet}) \to B(C_1'[1]^{\bullet}) \otimes \widetilde{D}'[1]^{\bullet} \otimes B(C_0'[1]^{\bullet}) 
\]
as the homomorphism determined by 
$\{{\mathfrak f}_k^{(i)}\}$, $i=0,1$, and $\{{\phi}_{k_1,k_0}\}$.  
We call $\phi=\{{\phi}_{k_1,k_0}\}$ a weakly filtered 
$A_{\infty}$-homomorphism of filtered $A_{\infty}$-bimodules, if 
\[
\widehat{d}_{{\mathfrak n}'} \circ \widehat{\phi} = 
\widehat{\phi} \circ \widehat{d}_{\mathfrak n}.  
\]
When we can take $c=0$ , $\phi=\{{\phi}_{k_1,k_0}\}$ is called 
a filtered $A_{\infty}$-homomorphism.  
Suppose that $C_i^{\bullet},C_i^{' \bullet}$ are $G$-gapped.  
Let $G' \subset \R \times 2\Z$ be a $G$-set such that 
${\rm pr}_1|_G^{-1} ((-\infty,\lambda])$ is finite 
for any $\lambda \in \R$ 
and ${\rm pr}_1(G)$ is bounded from below.   
We say that $\phi=\{{\phi}_{k_1,k_0}\}$ is $G'$-gapped, if 
\[
\phi_{k_1,k_0}=\sum_{\beta' \in G'} T^{{\rm pr}_1(\beta')} 
e^{{\rm pr}_2(\beta')/2} \phi_{k_1,k_0,\beta'},
\]
where 
\[
\phi_{k_1,k_0,\beta'}:B_{k_1}(\overline{C}_1[1]^{\bullet}) \otimes 
\overline{D}[1]^{\bullet} \otimes B_{k_0}(\overline{C}_0[1]^{\bullet}) 
\to \overline{D}'[1]^{\bullet}.
\]

The homotopy theory between filtered $A_{\infty}$-homomorphisms of 
filtered $A_{\infty}$-bimodules is also developed in \cite{FOOO}.  

We also proved the Whitehead theorem for (filtered) 
$A_{\infty}$-bimodules.  

\begin{thm}[Theorem 21.35 \cite{FOOO}]
Let $\phi:D^{\bullet} \to D^{' \bullet}$ be a gapped 
filtered $A_{\infty}$-bimodule 
homomorphism over $({\mathfrak f}^{(1)},{\mathfrak f}^{(0)})$, where 
${\mathfrak f}^{(i)}:C_i^{\bullet} \to C_i^{' \bullet}$ 
are homotopy equivalences.  
Suppose that $\phi_{(0,0,{\mathbf 0})}$ is a chain homotopy equivalence.  
Then $\phi$ is a homotopy equivalence of filtered 
$A_{\infty}$-bimodules.
\end{thm}

\begin{rem}
Here we require $\phi$ is a filtered $A_{\infty}$-homomorphism.  
Since a weakly filtered $A_{\infty}$-homomorphism, 
$\phi_{(0,0,{\mathbf 0})}$ may not induce a chain map with respect to 
${\mathfrak n}_{0,0,{\mathbf 0}}$ and 
${\mathfrak n}_{0,0,{\mathbf 0}}'$.
\end{rem}

\smallskip
\noindent
{\bf 2.6) Filtered $A_{n,K}$ structure.}

For a later argument in section 4, 
we recall the notion of filtered $A_{n,K}$-algebras.  
Let $G \subset \R{\geq 0} \times 2\Z$ be a monoid as above and 
$\beta_0={\mathbf 0} \in G$.    
For $\beta \in G$, we define 

\[
\parallel \beta \parallel = \left\{
\begin{array}{ll} 
\sup \{ n | \exists \beta_i \in G \setminus \{\beta_0\}, \ 
\sum_{i=1}^n \beta_i = \beta \} + [{\rm pr}_1(\beta)] -1 & 
\text{ if } \beta \neq \beta_0 \\
-1 & \text{ if } \beta = \beta_0,   
\end{array}
\right.
\]

Then we introduce a partial order on $(G \times \Z_{\geq 0})  
\setminus \{(\beta_0,0)\}$ by 
$(\beta_1,k_1) \succ (\beta_2,k_2)$ if and only if 
either 
\[
\parallel \beta_1 \parallel + k_1 > 
\parallel \beta_2 \parallel + k_2
\] or 
\[
\parallel \beta_1 \parallel + k_1 = 
\parallel \beta_2 \parallel + k_2, \text{ and } 
\parallel \beta_1 \parallel > \parallel \beta_2 \parallel.
\]
We write $(\beta_1,k_1) \sim (\beta_2,k_2)$, when 
\[
\parallel \beta_1 \parallel + k_1 = 
\parallel \beta_2 \parallel + k_2, \text{ and } 
\parallel \beta_1 \parallel = \parallel \beta_2 \parallel.
\]
We define $(\beta_1,k_1) \succsim (\beta_2,k_2)$ if 
either $(\beta_1,k_1) \succ (\beta_2,k_2)$ or $(\beta_1,k_1) \sim 
(\beta_2,k_2)$.  

We also write $(\beta,k) \prec (n',k')$, when 
$\parallel \beta \parallel + k < n' + k'$ or 
$\parallel \beta \parallel + k = n' + k'$ and 
$\parallel \beta \parallel < n'$.  

Let $\overline{C}^{\bullet}$ be a cochain complex over $R$ and 
$C^{\bullet}=\overline{C}^{\bullet}\otimes \Lambda_{0,nov}$.  
Suppose that there are 
\[
{\mathfrak m}_{k,\beta}:B(\overline{C}[1]^{\bullet}) \to \overline{C}[1]^{\bullet}
\]
for $(\beta,k) \in (G \times \Z) \setminus \{(\beta_0,0)\}$ 
with $(\beta,k) \prec (n,K)$.  
We also suppose that ${\mathfrak m}_{1,\beta_0}$ is 
the boundary operator of the cochain complex $C^{\bullet}$.  
We call $(C^{\bullet},\{{\mathfrak m}_{k,\beta}\})$ 
a $G$-gapped filtered $A_{n,K}$-algebra, if the following holds
\[
\sum_{\beta_1+\beta_2=\beta, k_1+k_2=k+1} \sum_i 
(-1)^{\deg'{\mathbf x}_i^{(1)}} {\mathfrak m}_{k_2,\beta_2}
\bigl({\mathbf x}_i^{(1)}, 
{\mathfrak m}_{k_1,\beta_1}({\mathbf x}_i^{(2)}), 
{\mathbf x}_i^{(3)} \bigr) = 0
\] 
for all $(\beta,k) \prec  (n,K)$, 
where 
\[
\Delta^2({\mathbf x})=\sum_i{\mathbf x}_i^{(1)} \otimes 
{\mathbf x}_i^{(2)} \otimes {\mathbf x}_i^{(3)}.
\]
Here $\Delta$ is the coproduct of the tensor coalgebra.  

We also have the notion of filtered $A_{n,K}$-homomorphisms, 
filtered $A_{n,K}$-homotopy equivalences in a natural way.  
In \cite{FOOO}, we proved the following:

\begin{thm}[Theorem 30.72 in \cite{FOOO}]\label{ext(n,K)}
Let $C_1^{\bullet}$ be a filtered $A_{n,K}$-algebra and 
$C_2^{\bullet}$ a filtered $A_{n',K'}$-algebra such that 
$(n,K) \prec (n',K')$.  
Let ${\mathfrak h}:C_1^{\bullet} \to C_2^{\bullet}$ 
be a filtered $A_{n,K}$-homomorphism.  
Suppose that ${\mathfrak h}$ is a filtered $A_{n,K}$-homotopy 
equivalence.  
Then there exist a filtered $A_{n',K'}$-algebra structure on 
$C_1^{\bullet}$ 
extending the given filtered $A_{n,K}$-algebra structure 
and a filtered $A_{n',K'}$-homotopy equivalence 
$C_1^{\bullet} \to C_2^{\bullet}$ extending the given 
filtered $A_{n,K}$-homotopy equivalence ${\mathfrak h}$.  
\end{thm}

\section{Canonical models}
In this section, we give the notion of canonical models and 
explain their construction after section 23, Chapter 5 of \cite{FOOO}.  
The unfiltered version of such a result goes back to Kadeishvili 
\cite{Kad}.  
There are two methods to construct canonical models in unfiltered 
case.  
One is based on obstruction theory due to Kadeishvili and 
the other uses the summation over trees due to Kontsevich and 
Soibelman \cite{KS}.  
Our argument is an adaptation of the latter argument and we also 
constructed the canonical models for filtered case.  

When the ground coefficent ring $R$ is a field, we have the 
following:  

\begin{thm}[Theorem 23.1, Theorem 23.2 in \cite{FOOO}]\label{canmodel}
(1) Any unfiltered $A_{\infty}$-algebra 
$(\overline{C}^{\bullet},\{\overline{\mathfrak m}_k\})$ is homotopy equivalent to an $A_{\infty}$-algebra 
$(\overline{C'}^{\bullet},\{\overline{\mathfrak m}'_k\})$ with 
$\overline{\mathfrak m}'_1=0$.  

\noindent
(2) Any gapped filtered $A_{\infty}$-algebra 
$(C^{\bullet},\{{\mathfrak m}_k\})$ 
is homotopy equivalent to a gapped filtered $A_{\infty}$-algebra 
$(C^{' \bullet},\{{\mathfrak m}'_k\})$ 
with $\overline{\mathfrak m}'_1=0$.  
Moreover, the homotopy equivalence can be taken as a gapped 
$A_{\infty}$-homomorphism.  
\end{thm}
 
An $A_{\infty}$-algebra is called {\it canonical}, if 
$\overline{\mathfrak m}_1=0$.  
A canonical model of an $A_{\infty}$-algebra is 
a canonical $A_{\infty}$-algebra homotopy equivalent to 
the original one.  
The statement (1) is Kadeishvili's theorem and implies that the 
$\overline{\mathfrak m}_1$-cohomology has a structure of an 
$A_{\infty}$-algebra.  
Note that, in general,  
we do not have ${\mathfrak m}_1$-cohomology, 
since ${\mathfrak m}_1 \circ {\mathfrak m}_1$ may not be zero.  
A filtered $A_{\infty}$-algebra is called {\it canonical}, if 
${\mathfrak m}_{1,0}=\overline{\mathfrak m}_1=0$.  
A canonical model of a filtered $A_{\infty}$-algebra is 
a canonical filtered $A_{\infty}$-algebra homotopy equivalent to 
the original one.  
 
Pick a submodule 
${\mathcal H}^{\bullet} \stackrel{\iota}{\hookrightarrow} 
\ker \overline{\mathfrak m}_1 \cap \overline{C}^{\bullet}$ such that 
$\iota_*: {\mathcal H}^k \cong 
{\rm H}^k(\overline{C}^{\bullet},\overline{\mathfrak m}_1)$, 
and $\Pi^k:\overline{C}^k \to {\mathcal H}^k \subset \overline{C}^k$ 
such that $\Pi^k \circ \Pi^k = \Pi^k$ and 
$\Pi^k \circ \overline{\mathfrak m}_1 = 0$.  
We will construct a structure of a filtered $A_{\infty}$-algebra on 
${\mathcal H}[1]^{\bullet} \otimes \Lambda_{0,nov}$ and 
a filtered $A_{\infty}$-homomorphism from 
${\mathcal H}[1]^{\bullet} \otimes \Lambda_{0,nov}$ 
to $C[1]^{\bullet}$, which is a weak homotopy equivalence.  
Since $R$ is a field, any cochain homomorphism inducing an isomorphism 
on cohomologies is a weak homotopy equivalence.  
Firstly, we observe the following:  

\begin{lem}\label{H}
There exist $G^k:\overline{C}^k \to \overline{C}^{k-1}$, 
$k=0,1,\dots, n$, such that 
\begin{align}
id - \Pi^k & =  -(\overline{\mathfrak m}_1 \circ G^k + G^{k+1} \circ \overline{\mathfrak m}_1), \label{htpy} \\
G^{k} \circ G^{k+1} & =  0. \label{sqzero} 
\end{align}
\end{lem}

From now on, let ${\mathcal H}^{\bullet} 
\stackrel{\iota}{\hookrightarrow} \overline{C}^{\bullet}$ 
be a subcomplex and $\Pi:\overline{C}^k \to 
\overline{C}^k$ be a projection to ${\mathcal H}^k$ such that 
there exist $G^k:\overline{C}^k \to \overline{C}^{k-1}$ 
satisfying (\ref{htpy}), (\ref{sqzero}) in Lemma \ref{H}.  
We do not assume that $\overline{\mathfrak m}_1|_{\mathcal H}=0$.  
Thus ${\mathcal H}^{\bullet}$ is not necessarily isomorphic 
to the cohomology $H^{\bullet}(\overline{C}^{\bullet})$.  
But the condition (\ref{htpy}) implies that 
$\iota_*: {\rm H}^{\bullet}({\mathcal H},\overline{\mathfrak m}_1|_{{\mathcal H}}) \cong 
{\rm H}^{\bullet}(\overline{C}^{\bullet},\overline{\mathfrak m}_1)$.  
Theorem \ref{canmodel} follows from the following:

\begin{thm}\label{reduction}
(1) There exists a structure 
$\{\overline{\mathfrak m}_k'\}_{k=1}^{\infty}$ of 
an $A_{\infty}$-algebra on ${\mathcal H}$ with 
${\mathfrak m}_1'={\mathfrak m}_1|_{\mathcal H}$.  
The inclusion $\iota$ extends to an $A_{\infty}$-homomorphism 
$\{\overline{\mathfrak f}_k\}_{k=1}^{\infty}$ 
with $\overline{\mathfrak f}_1 = \iota$.  

\noindent
(2) There exists a structure 
$\{{\mathfrak m}_k\}_{k=0}^{\infty}$ of a filtered $A_{\infty}$-algebra 
on ${\mathcal H} \otimes \Lambda_{0,nov}$.  
The inclusion $\iota$ extends to a filtered $A_{\infty}$-homomorphism 
$\{{\mathfrak f}\}_{k=0}^{\infty}$ with 
${\mathfrak f}_{1,0}=\iota$.  
\end{thm}

Let $G$ be a monoid as in section 2 and ${\rm pr}_1(G)=\{\lambda_{(i)}\}$ 
such that 
\[
0=\lambda_{(0)} < \lambda_{(1)} < \lambda_{(2)} < \cdots \to +\infty, 
\]
unless $G=\{(0,0)\}$.  
We write 
\[
{\mathfrak m}_{k,i}=
\sum_{\beta \in G \ \ {\rm pr}_1(\beta)=\lambda_{(i)}} 
e^{{\rm pr}_2(\beta)/2} {\mathfrak m}_{k,\beta}
\]
and
\[
{\mathfrak m}_{k,i}^{\circ}=T^{\lambda_{(i)}}{\mathfrak m}_{k,i}.
\]
Thus ${\mathfrak m}_k= \sum_i {\mathfrak m}_{k,i}^{\circ}$. 
Here ${\mathfrak m}_{k,i}$ is considered as 
\[
{\mathfrak m}_{k,i}:B_k(\overline{C}[1]^{\bullet}) \otimes R[e,e^{-1}] 
\to \overline{C}[1]^{\bullet} \otimes R[e,e^{-1}].
\]  
By an abuse of notation, we also denote by $\Pi^k$, $G^k$ 
the extensions thereof 
to $\overline{C}^k \otimes R[e,e^{-1}]$ as a $R[e,e^{-1}]$-module 
homomorphism.  

In order to define a $G$-gapped filtered $A_{\infty}$-structure 
on ${\mathcal H}[1]^{\bullet} \otimes \Lambda_{0,nov}$ and 
a $G$-gapped $A_{\infty}$-homomorphism from ${\mathcal H}[1]^{\bullet} 
\otimes \Lambda_{0,nov}$ to $C[1]^{\bullet}$, we introduce some notation.  

A decorated planar rooted tree is a quintet  
$\Gamma=(T,i,v_0,V_{tad}, \eta)$, which consists of 
\begin{itemize}
\item $T$ is a tree, 
\item $i:T \to D^2$ is an embedding, 
\item $v_0$ is the root vertex, 
\item $V_{tad}=\{ \text{vertices of valency } 1\} 
\setminus C^0_{ext}(T)$,
\item $\eta:C^0_{int}(T)=C^0(T) \setminus C^0_{ext}(T) \to \{0,1,2,\dots\}.$
\end{itemize}
Here $C^0(T)$ is the set of vertices of the tree $T$, 
$C^0_{ext}(T)=i^{-1}(\partial D^2)$ is the set of exterior vertices 
and $C^0_{int}(T)$ is the set of interior vertices.  
Note that the root vertex $v_0$ is an exterior vertex and 
$V_{tad} \subset C^0_{int}(T)$.  
Let ${G}_k^+$ be the set of $\Gamma=(T,i,v_0,V_{tad},\eta)$ 
such that $\# C^0_{ext} = k$ and $\eta(v) > 0$ if $v \in C^0_{int}(T)$ 
is a vertex of valency $1$ or $2$.  
We set $E(\Gamma)=\sum_{v \in C^0_{int}(T)} \lambda_{(\eta(v))}$.  

For each $\Gamma \in {G}_{k+1}^+$, we construct 
\[
{\mathfrak m}_{\Gamma}:B_k({\mathcal H}[1]^{\bullet}) \otimes R[e,e^{-1}] \to 
{\mathcal H}[1]^{\bullet} \otimes R[e,e^{-1}],
\] 
which is of degree $1$ and 
\[
{\mathfrak f}_{\Gamma}:B_k({\mathcal H}[1]^{\bullet}) \otimes R[e,e^{-1}] \to 
\overline{C}[1]^{\bullet} \otimes R[e,e^{-1}],
\]
which is of degree $0$.  
Then we define 
\[
{\mathfrak m}_k'=\sum_{\Gamma \in {G}_{k+1}^+} 
T^{E(\Gamma)} {\mathfrak m}_{\Gamma}:
B_k({\mathcal H}[1]^{\bullet}) \otimes \Lambda_{0,nov} \to 
{\mathcal H}[1]^{\bullet} \otimes \Lambda_{0,nov}
\]
and 
\[
{\mathfrak f}_k=\sum_{\Gamma \in {G}_{k+1}^+}
T^{E(\Gamma)} {\mathfrak f}_{\Gamma}:
B_k({\mathcal H}[1]^{\bullet}) \otimes \Lambda_{0,nov} \to 
\overline{C}[1]^{\bullet} \otimes \Lambda_{0,nov}.
\]
We will show that 
$({\mathcal H}[1]^{\bullet} \otimes \Lambda_{0,nov}, \{{\mathfrak m}_k'\})$ 
is a $G$-gapped filtered $A_{\infty}$-algebra and 
${\mathfrak f}=\{{\mathfrak f}_k\}$ is a $G$-gapped 
$A_{\infty}$-homomorphism, which is a weak homotopy equivalence.  
Then the Whitehead type theorem implies that ${\mathfrak f}$ is 
a homotopy equivalence.    

\smallskip \noindent
{\bf Step 1.} \ \ The case that $\# C^0_{int}(T)=0$.

Such a $T$ consists of two exterior vertices and an edge 
joining them.  
Therefore, there is unique element $\Gamma_0$, which belongs to 
$G_2^+$.  
We define 
\[
{\mathfrak m}_{\Gamma_0}=
\overline{\mathfrak m}_1|_{{\mathcal H}[1]^{\bullet}}
\] 
and 
\[
{\mathfrak f}_{\Gamma_0}:{\mathcal H}[1]^{\bullet} \otimes R[e,e^{-1}] \to 
\overline{C}[1]^{\bullet} \otimes R[e,e^{-1}]
\]
to be the inclusion $\iota$.  

\smallskip \noindent
{\bf Step 2.} \ \ The case that $\# C^0_{int}(T)=1$.

\begin{figure}
\begin{center}
\include{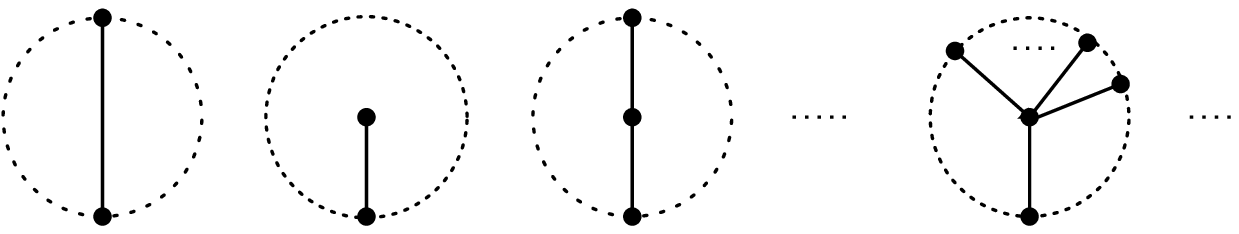}
\end{center}
\caption{}
\label{elmdectree}
\end{figure}

For any $k=0,1,2, \dots$, 
there is a unique planar tree with $\#C^0_{ext}(T)=k+1$ and 
$\#C^0_{int}(T)=1$.   
Let $\Gamma_{k+1} \in G_{k+1}^+$ be a decorated planar tree with 
one interior vertex $v$, see Figure \ref{elmdectree}.  

We define 
\[
{\mathfrak m}_{\Gamma_{k+1}}=\Pi \circ {\mathfrak m}_{k,\eta(v)}:
B_k({\mathcal H}[1]^{\bullet}) \otimes R[e,e^{-1}] \to 
{\mathcal H}[1]^{\bullet} \otimes R[e,e^{-1}]
\]
and 
\[
{\mathfrak f}_{\Gamma_{k+1}}=G \circ {\mathfrak m}_{k,\eta(v)}: 
B_k({\mathcal H}[1]^{\bullet}) \otimes R[e,e^{-1}] \to 
\overline{C}[1]^{\bullet} \otimes R[e,e^{-1}].
\]
Since the degree of $\Pi$, resp. $G$, is $0$, resp. $-1$, 
${\mathfrak m}_{\Gamma_{k+1}}$, resp. ${\mathfrak f}_{\Gamma_{k+1}}$, 
is of degree $1$, resp. $0$.  

\smallskip \noindent
{\bf Step 3.} \ \ General case. 

\begin{figure}
\begin{center}
\include{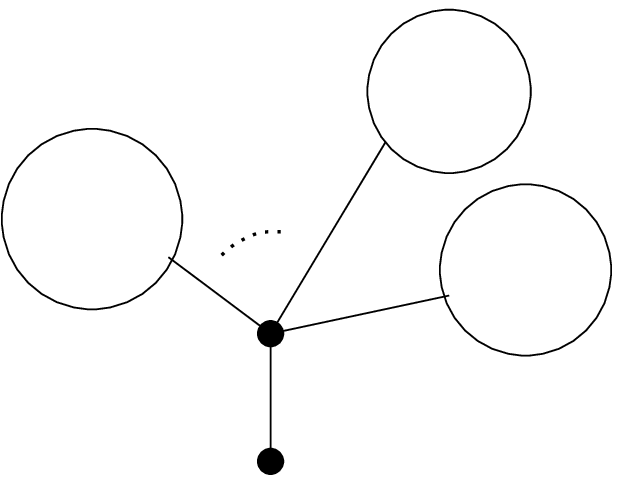}
\end{center}
\caption{}
\label{gendectree}
\end{figure}

Let $v_1$ is the vertex closest to the root vertex $v_0$.  
Cut the decorated planar tree at $v_1$, then $\Gamma$ is decomposed 
into decorated planar trees $\Gamma^{(1)}, \dots, \Gamma^{({\ell})}$ 
and an interval toward $v_0$ in counter-clockwise order, see 
Figure \ref{gendectree}.  

Then we define 
\[
{\mathfrak m}_{\Gamma}=\Pi \circ {\mathfrak m}_{\ell, \eta(v_1)} 
\circ ({\mathfrak f}_{\Gamma^{(1)}} \otimes \cdots \otimes {\mathfrak f}_{\Gamma^{({\ell})}})
\]
and 
\[
{\mathfrak f}_{\Gamma}=G \circ {\mathfrak m}_{\ell, \eta(v_1)} 
\circ ({\mathfrak f}_{\Gamma^{(1)}} \otimes \cdots \otimes {\mathfrak f}_{\Gamma^{({\ell})}}).
\]

\smallskip

Finally we define 
\[
{\mathfrak m}'_k=\sum_{\Gamma \in G_{k+1}^+} 
T^{E(\Gamma)} {\mathfrak m}_{\Gamma}
\] 
and 
\[
{\mathfrak f}_k=\sum_{\Gamma \in G_{k+1}^+} 
T^{E(\Gamma)} {\mathfrak f}_{\Gamma}.
\] 
As in \S 2.1, we obtain a graded coderivation 
\[
\widehat{d}'=\sum_k \widehat{\mathfrak m}'_k:B({\mathcal H}[1]^{\bullet}) \otimes \Lambda_{0,nov} \to B({\mathcal H}[1]^{\bullet}) \otimes \Lambda_{0,nov}
\] 
and a (formal) coalgebra homomorphism 
\[
\widehat{\mathfrak f}:B({\mathcal H}[1]^{\bullet}) \otimes \Lambda_{0,nov} \to 
B(C[1]^{\bullet}).
\] 

We will show the following:

\begin{prop}\label{mainprop}
\[
\widehat{\mathfrak f} \circ \widehat{d}' = \widehat{d} \circ 
\widehat{\mathfrak f},
\]
where $\widehat{d}=\sum_k \widehat{\mathfrak m}_k:B(C[1]^{\bullet}) \to B(C[1]^{\bullet})$. 
\end{prop}

Since $\overline{\mathfrak f}_1={\mathfrak f}_{\Gamma_0}$ is the 
inclusion, we find that $\widehat{\mathfrak f}$ is injective 
using the energy filtration and the number filtration on the bar 
complex.  
Then $\widehat{d}' \circ \widehat{d}' = 0$ follows from 
$\widehat{d} \circ \widehat{d} = 0$.  
Hence we obtain the following:

\begin{cor}
(1) $({\mathcal H}[1]^{\bullet} \otimes \Lambda_{0,nov}, \{{\mathfrak m}'_k\})$ 
is a $G$-gapped filtered $A_{\infty}$-algebra.  

\noindent
(2) $\widehat{\mathfrak f}$ is a $G$-gapped $A_{\infty}$-homomorphism 
from $({\mathcal H}[1]^{\bullet} \otimes \Lambda_{0,nov}, \{{\mathfrak m}'_k\})$ 
to $(C[1]^{\bullet}, \{{\mathfrak m}_k\})$.  
\end{cor}

The rest of this section is devoted to the proof of 
Proposition \ref{mainprop}, 
which is equivalent to that 
\[
{\mathfrak f} \circ \widehat{d}' = {\mathfrak m} \circ 
\widehat{\mathfrak f}  
\]
as maps $B({\mathcal H}[1]^{\bullet}) \otimes \Lambda_{0,nov} \to C[1]^{\bullet}$, 
where 
\[
{\mathfrak f}:B({\mathcal H}[1]^{\bullet}) \otimes \Lambda_{0,nov} 
\stackrel{\widehat{\mathfrak f}}{\to} B(C[1]^{\bullet}) \stackrel{\rm pr}{\to} 
C[1]^{\bullet},
\]
and 
\[
{\mathfrak m}:B(C[1]^{\bullet}) \stackrel{\widehat{d}}{\to} B(C[1]^{\bullet}) 
\stackrel{\rm pr}{\to} C[1]^{\bullet}.
\]
Namely, ${\mathfrak f}\vert_{B_k({\mathcal H}[1]^{\bullet}) \otimes \Lambda_{0,nov}} 
= {\mathfrak f}_k$, 
${\mathfrak m}\vert_{B_k({\mathcal H}[1]^{\bullet}) \otimes \Lambda_{0,nov}} 
= {\mathfrak m}_k$.  

We introduce an order on $\{(k,i) \vert k,i =0,1,2, \dots \}$ by 
$(k_1,i_1) \prec (k_2,i_2)$ if 
either $i_1 < i_2$ or $i_1=i_2$ and $k_1 < k_2$.  
We show the following claim by the induction on $(k,i)$.  

\begin{claim}
\[
{\mathfrak f} \circ \widehat{d}' \equiv {\mathfrak m} \circ 
\widehat{\mathfrak f}  \mod T^{\lambda_{(i+1)}}\cdot C[1]^{\bullet} \text{ \rm on } 
B_k({\mathcal H}[1]^{\bullet}) \otimes \Lambda_{0,nov}.
\]
\end{claim}

\begin{figure}
\begin{center}
\include{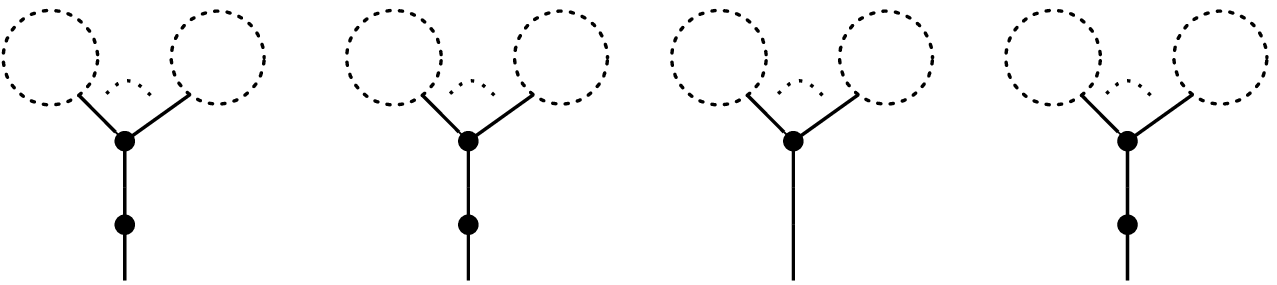}
\end{center}
\caption{}
\label{green}

\begin{center}
\include{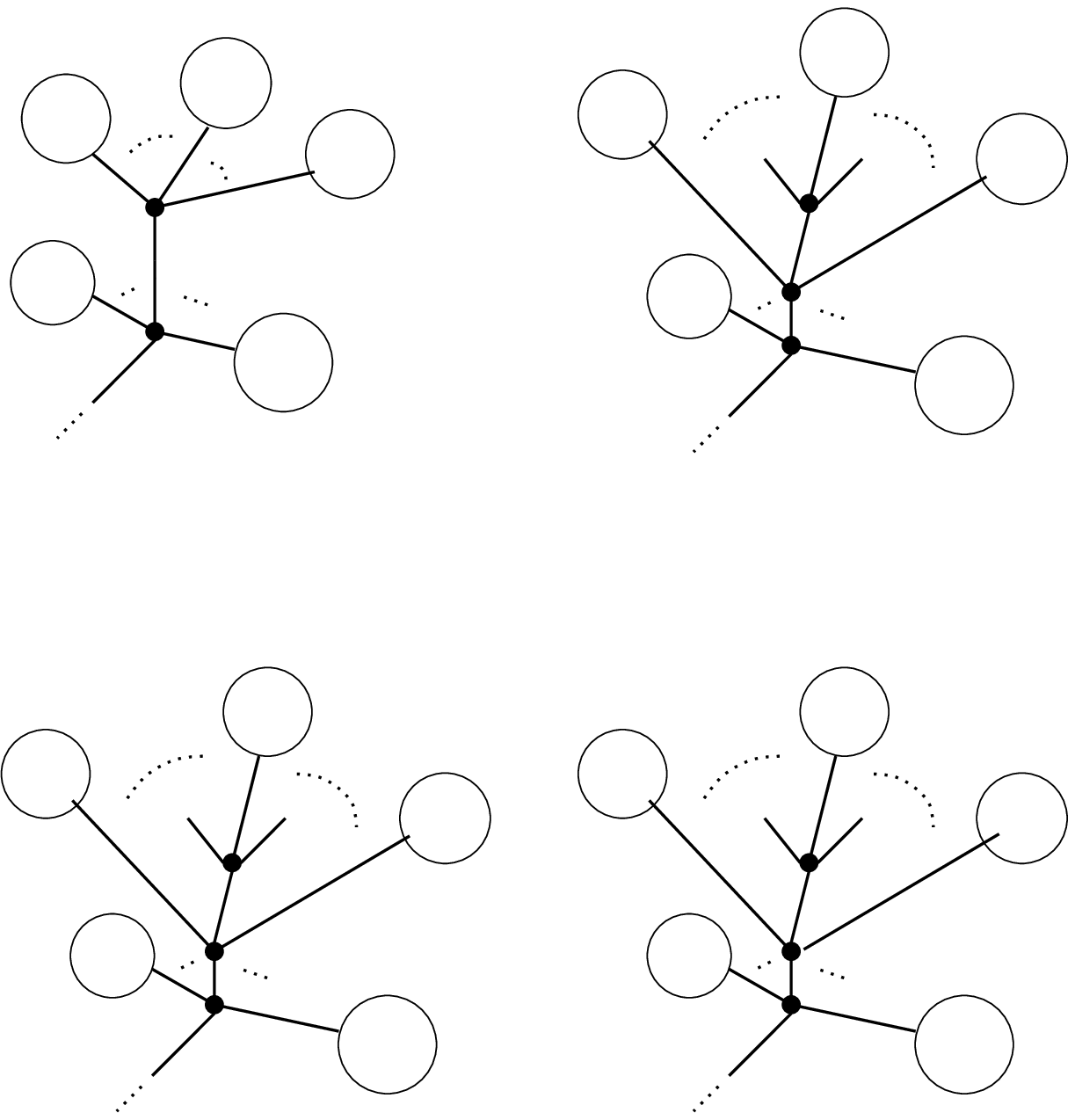}
\end{center}
\caption{}
\label{ainfty0}
\end{figure}

The key ingredients in the proof are the following relations 
presented in Figures \ref{green}, \ref{ainfty0}. 

Firstly, we consider the case that $i=0$.  
Claim $(0,0)$ follows from the gapped condition.  
By the choice of ${\mathcal H}$, Claim (1,0) holds clearly.  
Suppose that Claim $(\ell,0)$ holds for $\ell < k$.  
Note that 
\[
\overline{\mathfrak f}_k={\mathfrak f}_{k,0}^{\circ}  
= \bigl(\sum_{1<\ell\leq k} G \circ \overline{\mathfrak m}_{\ell} \circ 
\widehat{\overline{\mathfrak f}} + 
\delta_{k1} {\mathfrak f}_{\Gamma_0}\bigr)|_{B_k({\mathcal H}[1]^{\bullet})}
\] 
and 
\[
\overline{\mathfrak m}'_k={\mathfrak m}_{k,0}^{\circ}
= \bigl(\sum_{1<\ell\leq k} \Pi \circ \overline{\mathfrak m}_{\ell} 
\circ \widehat{\overline{\mathfrak f}} + 
\delta_{k1} {\mathfrak m}_{\Gamma_0}\bigr)|_{B_k({\mathcal H}[1]^{\bullet})}.
\]
Here $\delta_{ij}$ is Kronecker's delta.  
Recall that ${\mathfrak m}_{\Gamma_0}={\mathfrak m}_{1,0}|_{\mathcal H}$ and 
${\mathfrak f}_{\Gamma_0}$ is the inclusion.  
Note also that the restriction of $\widehat{\overline{\mathfrak f}}$ 
to $B_k({\mathcal H}[1]^{\bullet})$ in the right hand sides is determined by 
$\overline{\mathfrak f}_1, \dots, \overline{\mathfrak f}_{k-1}$.  

Thus we have 

\begin{align}
\overline{\mathfrak m} \circ \widehat{\overline{\mathfrak f}}
|_{B_k({\mathcal H}[1]^{\bullet})} 
& = 
\bigl(\overline{\mathfrak m}_1 \circ \widehat{\overline{\mathfrak f}} + 
\sum_{1< \ell \leq k} \overline{\mathfrak m}_{\ell} \circ 
\widehat{\overline{\mathfrak f}}\bigr)|_{B_k({\mathcal H}[1]^{\bullet})}  \nonumber \\
& =  
\bigl(\sum_{1< \ell \leq k} \overline{\mathfrak m}_1 \circ G \circ 
\overline{\mathfrak m}_{\ell} \circ \widehat{\overline{\mathfrak f}} + 
\delta_{k1} \overline{\mathfrak m}_1 \circ {\mathfrak f}_{\Gamma_0} 
+ \sum_{1< \ell \leq k} \overline{\mathfrak m}_{\ell} 
\circ \widehat{\overline{\mathfrak f}}\bigr)|_{B_k({\mathcal H}[1]^{\bullet})}.  
\nonumber \\
& = 
\bigl(\sum_{1< \ell \leq k} \Pi \circ \overline{\mathfrak m}_{\ell} 
\circ \widehat{\overline{\mathfrak f}} - 
\sum_{1< \ell \leq k} G \circ \overline{\mathfrak m}_1 \circ 
\overline{\mathfrak m}_{\ell} \circ \widehat{\overline{\mathfrak f}}
+\delta_{k1} \overline{\mathfrak m}_1 \circ {\mathfrak f}_{\Gamma_0} 
\bigr)
|_{B_k({\mathcal H}[1]^{\bullet})}
\nonumber \\
& =  
\bigl({\mathfrak f}_{\Gamma_0} \circ \overline{\mathfrak m}'_k - 
\sum_{1< \ell \leq k} G \circ \overline{\mathfrak m}_1 \circ 
\overline{\mathfrak m}_{\ell} \circ \widehat{\overline{\mathfrak f}}
\bigr)
|_{B_k({\mathcal H}[1]^{\bullet})} 
\nonumber \\
& = 
\bigl({\mathfrak f}_{\Gamma_0} \circ \overline{\mathfrak m}'_k + 
\sum_{1< \ell' \leq k} G \circ \overline{\mathfrak m}_{\ell'} \circ 
\widehat{\overline{d}} \circ \widehat{\overline{\mathfrak f}}\bigr)
|_{B_k({\mathcal H}[1]^{\bullet})}. 
\nonumber
\end{align}

Here we used the fact that 
$\overline{\mathfrak m}_1 \circ G + G \circ 
\overline{\mathfrak m}_1 = \Pi - id$ and the $A_{\infty}$-relation 
$\widehat{\overline{d}} \circ \widehat{\overline{d}} = 0$.   
Since we assumed Claim $(\ell,0)$ for $\ell < k$, i.e., 
\[
\overline{\mathfrak m} \circ \widehat{\overline{\mathfrak f}} = 
\overline{\mathfrak f} \circ \widehat{\overline{d}'} 
\text{ on } B_{\ell}({\mathcal H}[1]^{\bullet}),
\]  
we have 
\[
\widehat{\overline{d}} \circ \widehat{\overline{\mathfrak f}} \equiv 
\widehat{\overline{\mathfrak f}} \circ \widehat{\overline{d}'} 
\mod \overline{C}[1]^{\bullet}=B_1(\overline{C}[1]^{\bullet})
\text{ on } B_k({\mathcal H}[1]^{\bullet}).
\]
Therefore we find that 
\[
(\sum_{1< \ell' \leq k} G \circ \overline{\mathfrak m}_{\ell'} \circ 
\widehat{\overline{d}} \circ \widehat{\overline{\mathfrak f}})
|_{B_k({\mathcal H}[1]^{\bullet})}
= 
(\sum_{1< \ell' \leq k} G \circ \overline{\mathfrak m}_{\ell'} \circ 
\widehat{\overline{\mathfrak f}} \circ \widehat{\overline{d}'})
|_{B_k({\mathcal H}[1]^{\bullet})}.  
\]
Hence we showed Claim $(k,0)$, i.e.,   
\[
\overline{\mathfrak m} \circ \widehat{\overline{\mathfrak f}} 
= {\overline{\mathfrak f}} \circ \widehat{\overline{d}'}
\] 
on $B_k({\mathcal H}[1]^{\bullet})$.  

Next we assume that Claim $(k,i)$ holds for all $k=0,1,2, \dots$.  
We prove Claim $(k,i+1)$ by the induction on $k$.  
Note that Case 3-1 below does not occur in the case that $k=0$.  

First of all, we recall from the definition of 
$G_{k+1}^+$ that 
\[
{\mathfrak f}_k=\sum_{\Gamma \in G_{k+1}^+} T^{E(\Gamma)} 
{\mathfrak f}_{\Gamma} 
= \sum_{(\ell, j) \neq (1,0)} G \circ {\mathfrak m}_{\ell,j}^{\circ} 
\circ \widehat{\mathfrak f} |_{B_k({\mathcal H}[1]^{\bullet}) \otimes \Lambda_{0,nov}} + \delta_{k1} 
{\mathfrak f}_{\Gamma_0}.
\]
Then we have 

\begin{align}
{\mathfrak m} \circ \widehat{\mathfrak f}
|_{B_k({\mathcal H}[1]^{\bullet})} 
= &
\bigl({\mathfrak m}_{1,0} \circ \widehat{\mathfrak f} + 
\sum_{(\ell, j) \neq (1,0)} {\mathfrak m}_{\ell,j}^{\circ} \circ 
\widehat{\mathfrak f}\bigr)|_{B_k({\mathcal H}[1]^{\bullet})}  \nonumber \\
= & 
\bigl(\sum_{(\ell,j) \neq (1,0)} {\mathfrak m}_{1,0} \circ G \circ 
{\mathfrak m}_{\ell,j}^{\circ} \circ \widehat{\mathfrak f} + 
\delta_{k1} {\mathfrak m}_{1,0} \circ {\mathfrak f}_{\Gamma_0} 
\nonumber \\
& 
+ \sum_{(\ell,j) \neq (1,0)} {\mathfrak m}_{\ell,j}^{\circ} 
\circ \widehat{\mathfrak f}\bigr)|_{B_k({\mathcal H}[1]^{\bullet})}.  
\nonumber \\
= & 
\bigl(\sum_{(\ell,j) \neq (1,0)} \Pi \circ 
{\mathfrak m}_{\ell,j}^{\circ} \circ \widehat{\mathfrak f} - 
\sum_{(\ell,j) \neq (1,0)} G \circ {\mathfrak m}_{1,0} \circ 
{\mathfrak m}_{\ell,j}^{\circ} \circ \widehat{\mathfrak f} 
\nonumber \\
& 
+\delta_{k1} {\mathfrak m}_{1,0} \circ {\mathfrak f}_{\Gamma_0} 
\bigr)|_{B_k({\mathcal H}[1]^{\bullet})}
\nonumber \\
= & 
\bigl({\mathfrak f}_{\Gamma_0} \circ {\mathfrak m}'_k - 
\sum_{(\ell,j) \neq (1,0)} G \circ {\mathfrak m}_{1,0} \circ 
{\mathfrak m}_{\ell,j}^{\circ} \circ \widehat{\mathfrak f}
\bigr)
|_{B_k({\mathcal H}[1]^{\bullet})} 
\nonumber \\
= & 
\bigl({\mathfrak f}_{\Gamma_0} \circ {\mathfrak m}'_k + 
\sum_{(\ell',j') \neq (1,0)} G \circ 
{\mathfrak m}_{\ell',j'}^{\circ} 
\circ \widehat{d} \circ \widehat{\mathfrak f}\bigr)
|_{B_k({\mathcal H}[1]^{\bullet})}. 
\nonumber
\end{align}
In the third equality, we used the fact that 
${\mathfrak m}_{1,0} \circ G + G \circ 
{\mathfrak m}_{1,0} = \Pi - id$.   

We will show that 
\[
\sum_{(\ell',j') \neq (1,0)} G \circ {\mathfrak m}_{\ell',j'}^{\circ} 
\circ \widehat{d} \circ \widehat{\mathfrak f} \equiv 
\sum_{(\ell',j') \neq (1,0)} G \circ {\mathfrak m}_{\ell',j'}^{\circ} 
\circ \widehat{\mathfrak f} \circ \widehat{d'} \mod 
T^{\lambda_{(i+2)}},
\]
which implies that 
\[
{\mathfrak m} \circ \widehat{\mathfrak f} \equiv 
\widehat{\mathfrak f} \circ \widehat{d'} \mod T^{\lambda_{(i+2)}}.
\]

\smallskip \noindent
{\bf Case 1: $\ell'=0$.} \ \ 
Note that the $B_0(C[1]^{\bullet})=\Lambda_{0,nov}$-components of 
${\mathrm{Im}} \ \widehat{d} \circ \widehat{\mathfrak f}$ and 
${\mathrm{Im}} \ \widehat{\mathfrak f} \circ \widehat{d}'$ are zero.  
Hence we have
\[
{\mathfrak m}_{0,j'}^{\circ} 
\circ \widehat{d} \circ \widehat{\mathfrak f} =
{\mathfrak m}_{0,j'}^{\circ} 
\circ \widehat{\mathfrak f} \circ \widehat{d'} = 0.  
\]

\smallskip \noindent
{\bf Case 2: $\ell'=1$.} \ \ 
For $j' \neq 0$, ${\mathfrak m}_{1,j'}^{\circ} \equiv 0 
\mod T^{\lambda_{(1)}}$.  
By the induction hypothesis, we have 
\[
\widehat{\mathfrak f} \circ \widehat{d'} \equiv 
\widehat{d} \circ \widehat{\mathfrak f} \mod T^{\lambda_{i+1}}.
\]
Since $\lambda_{(i+2)} \leq \lambda_{(i+1)} + \lambda_{(1)}$, 
we obtain 
\[
{\mathfrak m}_{1,j'}^{\circ} \circ \widehat{\mathfrak f} \circ 
\widehat{d'} \equiv 
{\mathfrak m}_{1,j'}^{\circ} \circ \widehat{d} \circ 
\widehat{\mathfrak f} \mod T^{\lambda_{(i+2)}}.
\] 

\smallskip \noindent
{\bf Case 3: $\ell' \geq 2$.} \ \ 
Let ${\mathbf x} \in B_k({\mathcal H}[1]^{\bullet}) \otimes \Lambda_{0,nov}$.  
Write 
\[
\Delta^{\ell'-1} {\mathbf x} = \sum_a {\mathbf x}_{1,a} \otimes 
\cdots \otimes {\mathbf x}_{\ell',a},
\]
where $\Delta$ is the coproduct and ${\mathbf x}_{i,a} \in 
B_{k_{i,a}}({\mathcal H}[1]^{\bullet}) \otimes \Lambda_{0,nov}$.  
Then we have 

\begin{align}
\widehat{\mathfrak f} \circ \widehat{d'} ({\mathbf x}) 
= & \sum_a \sum_j (-1)^{\deg' {\mathbf x}_{1,a}+\cdots + 
\deg'{\mathbf x}_{j-1,a}}
{\mathfrak f}_{k_{1,a}}({\mathbf x}_{1,a}) \otimes \cdots \otimes 
{\mathfrak f}_{k_{j,a}}(\widehat{d'}({\mathbf x}_{j,a})) \otimes 
\cdots \nonumber \\
& \hspace{0.3in} \cdots \otimes {\mathfrak f}_{k_{\ell',a}}({\mathbf x}_{\ell',a}). 
\nonumber
\end{align}

\smallskip \noindent
{\bf Case 3-1: $k_{j,a} < k$.} 
In this case, we have 
\[
{\mathfrak f}_{k_{j,a}}(\widehat{d'}({\mathbf x}_{j,a}) \equiv 
{\mathfrak m} \circ \widehat{\mathfrak f}({\mathbf x}_{j,a}) 
\mod T^{\lambda_{(i+2)}}
\]
by the induction hypothesis.  
Hence 

\begin{align}
& 
{\mathfrak f}_{k_{1,a}}({\mathbf x}_{1,a}) \otimes \cdots \otimes 
{\mathfrak f}_{k_{j,a}}(\widehat{d'}({\mathbf x}_{j,a})) \otimes 
\cdots \otimes {\mathfrak f}_{k_{\ell',a}}({\mathbf x}_{\ell',a}) 
\nonumber \\ 
\equiv & 
{\mathfrak f}_{k_{1,a}}({\mathbf x}_{1,a}) \otimes \cdots \otimes 
{\mathfrak m} \circ \widehat{\mathfrak f}({\mathbf x}_{j,a}) \otimes 
\cdots \otimes {\mathfrak f}_{k_{\ell',a}}({\mathbf x}_{\ell',a}) 
 \ 
\mod T^{\lambda_{(i+2)}}. \nonumber 
\end{align}

\smallskip \noindent
{\bf Case 3-2: $k_{j,a} = k$.} 
In this case, $k_{j',a} = 0$ for $j' \neq j$, i.e., 
${\mathbf x}_{j',a} \in B_0({\mathcal H}[1]^{\bullet}) \otimes \Lambda_{0,nov}$.  
Without loss of generality, we may assume that ${\mathbf x}_{j',a}=1$ 
for $j' \neq j$.  

By the induction hypothesis, we have 
\[
{\mathfrak f}(\widehat{d'}({\mathbf x}_{j,a})) \equiv 
{\mathfrak m}(\widehat{\mathfrak f}({\mathbf x}_{j,a})) 
\mod T^{\lambda_{(i+1)}},
\]
which implies that 
\[
{\mathfrak f}_0(1) \otimes \cdot \otimes {\mathfrak f}
(\widehat{d'}({\mathbf x}_{j,a})) \otimes \cdots \otimes 
{\mathfrak f}_1(1)  \equiv  
{\mathfrak f}_0(1) \otimes \cdot \otimes {\mathfrak m} 
(\widehat{\mathfrak f}({\mathbf x}_{j,a})) \otimes \cdots \otimes 
{\mathfrak f}_1(1) 
\mod T^{\lambda_{(i+2)}}.
\]
Here we used ${\mathfrak f}_0(1) \equiv 0 \mod T^{\lambda_{(1)}}$. 

In sum, we obtain Claim $(k,i+1)$ for all $k$.  

By the construction, $\overline{\mathfrak f}_1$ is a chain homotopy 
equivalence ($\Pi$ is a homotopy inverse).  
Therefore, Theorem \ref{whitehead} implies that 
$\{{\mathfrak f}_k\}$ is a homotopy equivalence of 
filtered $A_{\infty}$-algebras.  

\section{Filtered $A_{\infty}$-algebra associated to Lagrangian 
submanifolds}

Let $(M,\omega)$ be a closed symplectic manifold and 
$L$ a Lagrangian submanifold.  
We only consider the case that $L$ is an embedded compact 
Lagrangian submanifold without boundary equipped with a 
relative spin structure, see \S 44 in \cite{FOOO}.   
We constructed a filtered $A_{\infty}$-algebra 
associated to $L$ in $(M,\omega)$.  
As we explained in section 2, the framework of (filtered) 
$A_{\infty}$-algebras, bimodules, etc. is adequate to formulate 
the condition under which Floer complex is obtained.  

In this section, we briefly recall the way of constructing 
filtered $A_{\infty}$-algebra associated to $L$.  
Although the readers may find Proposition 4.1 below too technical, 
we present it precisely so that we can explain how to modify it 
for the purpose of section 5.  

A naive idea of the construction is to use the moduli space 
of pseudo-holomorphic discs to {\it deform} 
the intersection products of chains in $L$ 
in a similar way to the quantum cohomology, where the intersection 
product on (co)homology is deformed by the moduli space 
of pseudo-holomorphic spheres, more precisely, stable maps of genus 0.  
Here appears a difference: while the moduli spaces of stable maps of 
genus 0 are (virtual) cycles, the moduli spaces of stable 
bordered stable maps are, in general, not (virtual) cycles, 
but with codimension 1 boundary (in the sense of Kuranishi structure).  
Therefore, we cannot restrict ourselves to cycles and forced to 
work with chains.  
However, the intersection product is not defined 
in chain level, e.g., the self intersection of chains.  
We start with a subcomplex of the singular 
chain complex such that the inclusion induces an isomorphism 
on homology.   
Then take {\it perturbed} intersection product 
of generators of the subcomplex and add them to get a larger 
subcomplex such that the inclusion induces an isomorphism on homology.  
Once we get such nested subcomplexes, we apply the argument 
in the proof of Theorem 3.3 to define the operation
$\overline{\mathfrak m}_2$ on a fixed subcomplex.
This multiplicative structure is not associative, but associative 
up to homotopy.  
So we proceed to constructed other operations 
$\overline{\mathfrak m}_k$ in a similar way, see 
Corollary 30.89 in section 30.6, \cite{FOOO} for a detailed argument.  
In this way, we obtain an $A_{\infty}$-algebra.  

For the construction of the filtered $A_{\infty}$-algebra, 
we include the effect from the moduli space of bordered stable 
maps.  
We need to take perturbations of the moduli spaces to define 
the operations not only perturbation in the intersection product 
mentioned above.  
Our strategy is to construct an $A_{n,K}$-algebra on 
$C_{(g)}(L)$, which is generated by $\chi_{(g)}$ in 
Proposition \ref{transv}, for a sufficiently large $g$.  
Then we use the obstruction theory to extend a filtered 
$A_{n,K}$-structure to a filtered $A_{n',K'}$-structure 
($(n,K) \precsim (n',K')$).  
The resulting filtered $A_{\infty}$-structure is unique up to 
homotopy, see \S 30 in Chapter 7, \cite{FOOO}.  

Let $\mu_L \in {\rm H}^2(M,L;\Z)$ be the Maslov class of the 
Lagrangian submanifold $L$.  
We introduce an equivalence relation $\sim$ on ${\rm H}_2(M,L;\Z)$ 
by $\beta_1 \sim \beta_2$ if and only if 
$\omega(\beta_1) = \omega(\beta_2)$ and 
$\mu_L(\beta_1)=\mu_L(\beta_2)$.  

Pick an almost complex structure $J$ compatible with $\omega$.  
Denote by ${\mathcal M}(\beta;L,J)$ the moduli space 
of bordered stable maps $u:(\Sigma, \partial \Sigma) \to (M,L)$ 
of genus $0$ representing $\beta$ and 
by ${\mathcal M}_{k+1}(\beta;L,J)$ be the moduli space of 
bordered stable maps in the class $\beta$ of genus $0$ 
with $k+1$ marked points 
$z_0,z_1, \dots, z_k$ on the regular part of $\partial \Sigma$.  
Denote by ${\mathcal M}_{k+1}^{\rm main}(\beta;L,J)$ the component, 
on which the marked points $z_0, z_1, \dots, z_k$ respect 
the {\it counter-clockwise} cyclic order on the boundary of bordered 
semi-stable curve of genus 0 with connected boundary.  
Let ${\mathfrak G}(L)$ be 
the monoid contained in $\Pi(M,L)$ generated by 
$\beta$ with ${\mathcal M}(\beta;L,J) \neq \emptyset$.  
We write $\beta_0=0 \in {\mathfrak G}(L)$.  

Our basic idea is as follows.  
For singular simplices $P_1, \dots, P_k$ in $L$, 
we consider the fiber product in the sense of Kuranishi structure 
\[
{\mathcal M}_{k+1}^{\rm main}(\beta;P_1, \dots, P_k) ={\mathcal M}_{k+1}
^{\rm main}(\beta;L,J)_{\mathbf{ev}} \times_{L^k} (P_1 \times \dots \times P_k), 
\]
where ${\mathbf{ev}}=(ev_1, \dots , ev_k)$ is the evaluation map at 
$z_1, \dots , z_k$.  
(For the orientation issue, see Chapter 9 \cite{FOOO}.)  
Then we would like to define 
\[
{\mathfrak m}_{k,\beta}(P_1, \dots, P_k) = (ev_0: 
{\mathcal M}_{k+1}^{\rm main}(\beta;P_1, \dots, P_k) \to L),
\]
where $ev_0$ is the evaluation at $z_0$.  

Note that ${\mathcal M}_{k+1}^{\rm main}(\beta)$ is not necessarily 
a manifold or an orbifold and that 
$ev_i$ are not necessarily submersions even if 
${\mathcal M}_{k+1}^{\rm main}(\beta)$ is such a nice space.  
In order to deal with this issue, we introduced the notion of 
Kuranishi structure \cite {FO}, see also Appendix in \cite{FOOO}.  
Here is a digression on Kuranishi structure.  

Let $X$ be a compact Hausdorff space.  
A Kuranishi structure on $X$ consists of 
a covering of $X$ by Kuranishi neighborhoods of the same 
{\it virtual} dimension and coordinate changes among them.  
A Kuranishi neighborhood around $p \in X$ is a quintet  
 $(V_p,E_p,\Gamma_p, s_p, \psi_p)$, where 
\begin{itemize}
\item
$V_p$ is a smooth manifold of finite dimension,
\item
$E_p$ is a real vector bundle over $V_p$ of finite rank, 
\item
$\Gamma_p$ is a finite group acting smoothly and effectively 
on $V_p$ and $E_p$ such that $E_p \to V_p$ is a $\Gamma_p$-equivariant 
vector bundle, 
\item
$s_p$ is a $\Gamma_p$-equivariant section of $E_p \to V_p$, 
\item 
$\psi_p$ is a homeomorphism from $s_p^{-1}(0)/\Gamma_p$ to 
a neighborhood of $p$ in $X$.
\end{itemize}
The vector bundle $E_p \to V_p$ is called the {\it obstruction} bundle 
and the section $s_p$ the Kuranishi map.  
We have coordinate changes among Kuranishi neighborhoods, see 
\cite{FO}, \cite{FOOO}.  
We require that $\dim V_p - {\rm rank}~ E_p$ does not depend on 
$p \in X$ and call it the {\it virtual} dimension of the space $X$ 
equipped with Kuranishi structure.  

The moduli spaces of stable maps, bordered stable maps carry 
Kuranishi structures, hence we can locally describe the moduli space 
as $s_p^{-1}(0)/\Gamma_p$ in the definition of Kuranishi neighborhoods. 
If $s_p$ is transversal to the zero section, the moduli space is 
locally an orbifold.  
In general, we cannot perturb $s_p$ to 
a $\Gamma_p$-equivariant section $s'_p$, 
which is transversal to the zero section.  
Instead of single valued sections,  
we consider perturbations by 
$\Gamma_p$-equivariant {\it multi-valued} sections, 
each branch of which is transversal to the zero section.  
Then we arrange them compatible under the coordinate change.  
In this way, we obtain {\it perturbed} moduli spaces.  

We take a multi-valued perturbation $\mathfrak s$ of 
Kuranishi maps for 
${\mathcal M}_{k+1}^{\rm main}(\beta;P_1, \dots, P_k)$ such that 
each branch of $\mathfrak s$ is transversal to the zero section.  
After taking a triangulation of 
the perturbed zero locus 
${\mathcal M}_{k+1}^{\rm main}(\beta;P_1, \dots, P_k)^{\mathfrak s}$ 
of $\mathfrak s$, we obtain a {\it virtual} chain 
$$
ev_0: {\mathcal M}_{k+1}^{\rm main}(\beta;P_1, \dots, P_k)^{\mathfrak s}
\to L.
$$  

To make this argument rigorous, we build a sequence of 
subcomplexes of the singular chain complex of $L$ and 
a series of operations ${\mathfrak m}_{k,\beta}^{geo}$.    
For details, see Chapter 7 in \cite{FOOO}.  
Here we briefly recall a part of it, in particular, 
the construction of a series of subcomplexes of singular 
chain complex of $L$.  
In section 5, we explain how to arrange this construction 
in relation with the Morse theory.  

In \S 30 in \cite{FOOO}, we constructed countable sets $\chi_g(L)$ 
of singular $C^{\infty}$-simplices on $L$.  
For a simplex $P \in \chi_g(L)$, we call $g$ the generation of $P$.   
Write 
\[
\chi_{(g)}=\bigcup_{g' \leq g} \chi_{g'}(L)
\] 
and denote by $C_{(g)}(L;R)$ the $R$-vector space generated by 
$\chi_{(g)}(L)$.  
Let $S(L;R)$ be the singular $C^{\infty}$-chain complex of $L$ 
with coefficients in $R$.

\smallskip
\noindent
{\bf Condition 1.} \ \  Any face of $P \in \chi_g(L)$ belongs to 
$\chi_{(g)}(L)$.  

\smallskip
\noindent
{\bf Condition 2.} \ \ The inclusion $C_{(g)}(L) \to S(L;R)$ 
induces an isomorphism on homology.  

\smallskip

For $\beta \in {\mathfrak G}(L)$, we define 

\[
\parallel \beta \parallel = 
\left\{ 
\begin{array}{ll}
\sup\{n|\exists \beta_1,\dots,\beta_n \in {\mathfrak G}(L) \setminus 
\{\beta_0\}, \sum_{i=1}^n \beta_i = \beta \} 
+[\omega(\beta)] -1 & 
\text{ if } \beta \neq \beta_0 \\
-1 & \text{ if } \beta = \beta_0
\end{array}
\right.
\]
Here $[\omega(\beta)]$ is the largest integer not greater than 
$\omega(\beta)$.  

By Gromov's compactness, 
the number of $\beta \in {\mathfrak G}(L)$ with 
$\parallel \beta \parallel \leq C$ is finite for any $C$.  

Next we introduce an additional data 
${\mathfrak d}:\{1, \dots ,k\} \to \Z_{\geq 0}$, which is called 
a decoration.  
For a pair $({\mathfrak d},\beta)$ such that 
${\mathcal M}_{k+1}^{\rm main}(\beta) \neq \emptyset$, we define 

\[
\parallel ({\mathfrak d},\beta) \parallel = 
\left\{ 
\begin{array}{ll} 
\max_{i \in \{1, \dots, k\}} {\mathfrak d}(i) + 
\parallel \beta \parallel + k & \text{ if } k \neq 0 \\
\parallel \beta \parallel & \text{ if } k = 0.  
\end{array}
\right.
\]

We will take the fiber product of 
${\mathcal M}_{k+1}^{\rm main}(\beta)$ and singular simplices 
$P_i$ in $L$.  
The decoration $\mathfrak d$ is introduced in order to 
include the generations of singular simplices $P_i$ into the data.  
When we emphasize that the decoration $\mathfrak d$ 
is equipped with the moduli space 
${\mathcal M}_{k+1}^{\rm main}(\beta)$, we denote it by 
${\mathcal M}_{k+1}^{{\rm main}, {\mathfrak d}}(\beta)$.  

\begin{prop}[Proposition 30.35 in \cite{FOOO}]\label{transv}
For any $\delta > 0$ and ${\mathcal K} > 0$, there exist 
$\chi_{(g)}(L)$, $g=0,\dots,{\mathcal K}$, and multisections 
${\mathfrak s}_{{\mathfrak d},k,\beta,\vec{P}}$ for 
$\parallel ({\mathfrak d},\beta) \parallel \leq {\mathcal K}$ 
with the following properties: 

\begin{itemize}
\item $\chi_{(g)}(L)$ satisfies Conditions 1 and 2 above.  
\item Let $P_i \in \chi_{{\mathfrak d}(i)}(L), i=1,\dots,k$.  
We put 
\[
{\mathcal M}_{k+1}^{{\rm main},{\mathfrak d}}(\beta;P_1,\dots,P_k)
= 
{\mathcal M}_{k+1}^{{\rm main},{\mathfrak d}}(\beta)\times_{L^k} 
\prod P_i
\]
and define a multisection 
${\mathfrak s}_{{\mathfrak d},k,\beta,\vec{P}}$ thereof.   
${\mathfrak s}_{{\mathfrak d},k,\beta,\vec{P}}$ is transversal to 
the zero section.  
\item If $g=\parallel ({\mathfrak d}, \beta) \parallel$, then 
\[
ev_{0*}\bigl({\mathcal M}_{k+1}^{{\rm main},{\mathfrak d}}
(\beta;P_1,\dots,P_k)^{{\mathfrak s}_{{\mathfrak d},k,\beta,\vec{P}}} 
\bigr)
\]
is decomposed into elements of $\chi_{(g)}(L)$.  
Here and henceforth we denote 
\[
{\mathcal M}_{k+1}^{{\rm main},{\mathfrak d}}
(\beta;P_1,\dots,P_k)^{{\mathfrak s}_{{\mathfrak d},k,\beta,\vec{P}}}
:= {\mathfrak s}_{{\mathfrak d},k,\beta,\vec{P}}^{-1}(0).
\]
\item The multisections ${\mathfrak s}_{{\mathfrak d},k,\beta,\vec{P}}$ 
satisfy certain compatibility conditions.
\item The zero locus ${\mathfrak s}_{{\mathfrak d},k,\beta,\vec{P}}^{-1}
(0)$ is in a $\delta$-neighborhood of the zero locus of the original 
Kuranishi map.  
\end{itemize}
\end{prop}

For the compatibility conditions in the above statement, see 
Conditions 30.38 and 30.44 in \cite{FOOO}.  

Now we explain the way of constructing the filtered 
$A_{\infty}$-algebra associated to $L$.  

We put 
\[
{\mathfrak m}_{k,\beta}^{geo}(P_1,\dots,P_k)
=(ev_0:{\mathcal M}_{k+1}^{{\rm main},{\mathfrak d}}(\beta;P_1,\dots,
P_k)^{{\mathfrak s}_{{\mathfrak d},k,\beta,\vec{P}}} \to L),
\]
when $P_i \in \chi_{({\mathfrak d}(i))}$, $i=1,\dots, k$.  
Then ${\mathfrak m}_{k,\beta}^{geo}(P_1,\dots,P_k)$ is decomposed 
into elements of $\chi_{(g)}$, 
where $g = \parallel ({\mathfrak d},\beta) \parallel$.  
Using the idea in section 3, we showed the following:  

\begin{prop}[Proposition 30.78 in \cite{FOOO}]\label{geo}
For any $g_0, n, K$, there exists $g_1 > g_0$ and a filtered 
$A_{n,K}$-structure ${\mathfrak m}_{k,\beta}$ on $C_{(g_1)}(L) 
\otimes \Lambda_{0,nov}$ such that 
\[
{\mathfrak m}_{k,\beta}(P_1,\dots,P_k)
={\mathfrak m}_{k,\beta}^{geo}(P_1,\dots,P_k),
\]
if $P_i \in \chi_{(g_0)}(L)$.  
\end{prop}

Combining Theorem \ref{ext(n,K)} and Proposition \ref{geo}, 
we can construct a filtered $A_{\infty}$-algebra associated to 
$L$, for details see \cite{FOOO}.  Hence we obtain 

\begin{thm}[Theorem 10.11 in \cite{FOOO}]
Let $L$ be a relatively spin Lagrangian submanifold.  
Then there exist a countably generated subcomplex $C(L)$ of 
the singular chain complex and a filtered $A_{\infty}$-algebra 
structure on $C(L)\otimes \Lambda_{0,nov}$.  
\end{thm}

We also proved that the homotopy type of 
the filtered $A_{\infty}$-algebra is unique.  

Applying the construction of canonical models in section 3, 
we obtain a filtered $A_{\infty}$-algebra structure on 
${\rm H}(L) \otimes \Lambda_{0,nov}$.  

Let $(L_0,L_1)$ be a relative spin pair of Lagrangian submanifolds.  
Assume that $L_0$ and $L_1$ intersect transversely.  
Then we have the following:

\begin{thm}
Let $D^{\bullet}$ be a free $\Lambda_{0,nov}$-module 
generated by $L_0 \cap L_1$.  
Then there exists a filtered $A_{\infty}$-bimodule structure 
over filtered $A_{\infty}$-algebras associated to $L_i$, $i=0,1$.  
\end{thm}

\section{Canonical models and Morse complexes}

In this section, we apply Theorem \ref{reduction} 
and reduce the filtered $A_{\infty}$-structure 
on $C^{\bullet}(L) \otimes \Lambda_{0,nov}$ to the Morse complex 
$CM^{\bullet}(f) \otimes \Lambda_{0,nov}$.

We pick a specific Morse function 
as follows.  
Choose and fix a triangulation ${\mathfrak T}$ of $L$.  
We may assume that the triangulation is sufficiently fine 
by taking subdivision.  
Pick a Morse function $f:L \to \R$ with the following property.  
Critical points of $f$ are in one-to-one correspondence with 
barycenters of simplices.  
Moreover, the Morse index of a critical point 
is equal to the dimension of the corresponding simplex.  
Then we can take a gradient-like vector field $X$ such that 
the unstable manifold $W^u(p)$ at each critical point $p$ is 
the interior of the corresponding simplex.  
Denote by $\{\rho_t\}$ the flow generated by $X$.  
(The function $f$ increases along the orbits of $\{\rho_t\}$.)  

Now we prove the following:

\begin{thm}\label{Morse}
Let $L$ be a relatively spin Lagrangian submanifold in a closed 
symplectic manifold $(M,\omega)$ and $f$ a Morse function on $L$ 
as above.  
Then Morse complex $CM^*(f)\otimes \Lambda_{0,nov}$ 
carries a structure of a filtered $A_{\infty}$-algebra, 
which is homotopy equivalent to the filtered $A_{\infty}$-algebra 
associated to $L$ constructed in \cite{FOOO}.  
\end{thm} 
The proof occupies the rest of this section.  
We explain how to choose $\chi_g(L)$ in section 4. 
Firstly, we choose and fix a linear order on the set of vertices 
in ${\mathfrak T}$.  
Then we regard each $T_i \in {\mathfrak T}$ as 
a singular simplex by the affine parametrization 
$\sigma_i:\Delta_{k_i} \to T_i$ 
respecting the order of the vertices.  
In particular, all simplices are oriented, hence 
the unstable manifolds $W^u(p)$.  
For our construction, we have to start with the following 
set of singular simplices.  
Set $\chi_{\mathfrak T}(L)=\{\sigma_i\}$ 
and identify the Morse complex $CM^{\bullet}(f)$ 
with $C_{\mathfrak T}(L)$, which is a subcomplex of the singular 
chain complex of $L$ generated by $\chi_{\mathfrak T}(L)$.  
Note that $\chi_{\mathfrak T}(L)$ satisfies Conditions 1 and 2 
given in section 4.  

We define $\chi_g(L) \supset \chi_{\mathfrak T}(L)$ 
in an inductive way as follows.  
For $g=-1$, we set $\chi_{-1}(L)=\chi_{\mathfrak T}(L)$.
For $g=0,1,\dots$, suppose that 
we constructed $\chi_{g'}(L)$, $g' < g$.  

We can choose the perturbations 
${\mathfrak s}_{{\mathfrak d},k,\beta,\vec{P}}$ in 
Proposition \ref{transv} with the following property.  

\smallskip

Each face $\tau$ of any simplex in the triangulation of 
\[
ev_{0*}\bigl({\mathcal M}_{k+1}^{{\rm main},{\mathfrak d}}
(\beta;P_1,\dots,P_k)^{{\mathfrak s}_{{\mathfrak d},k,\beta,\vec{P}}} 
\bigr)
\] 
with $g=\parallel ({\mathfrak d}, \beta) \parallel$
is transversal to the stable manifold $W^s(p)$ 
at any $p \in {\rm Crit}(f)$.  
Moreover, for each $\tau$ of dimension at most 
$\dim L$, there exists at most one 
$p=p(\tau) \in {\rm Crit}(f)$ such that the stable submanifold 
$W^s(p)$ is of complementary dimension to $\tau$ and 
$W^s(p)$ and $\tau$ intersect at a unique point.  
Denote by $T(p) \in {\mathfrak T}$ the simplex containing $p$.  
Let $\chi_g^{\circ}(L)$ be the set of these singular simplices 
$\tau$.

\smallskip

We have to add $\chi_g^{\circ}(L)$ to previous 
$\bigcup_{g' < q}\chi_{g'}(L)$.  
In order to guarantee Condition 2, we further add the following 
singular simplices to $\chi_g^{\circ}(L)$ and obtain 
$\chi_g(L)$.  
Denote by $\sigma^{\tau} \in \chi_{\mathfrak T}$ 
the singular simplex corresponding to $T(p(\tau))$.  
Define $\Pi(\tau)=\epsilon \sigma^{\tau}$, where 
$\epsilon = \pm 1$ is given by the following equation.  
\[
\tau \cap W^s(p(\tau)) = \epsilon W^u(p(\tau)) \cap W^s(p(\tau)),
\]
if there exists a unique stable manifold $W^s(p(\tau))$, 
which intersects $\tau$ transversely at a unique point. 
Otherwise, we define $\Pi(\tau)=0$.  
In particular, if $\tau > \dim L$, $\Pi(\tau)=0$.  
For each $\tau$ as above, we will find a singular chain $G(\tau)$ 
such that 
$$
\Pi (\tau) - \tau = \overline{\mathfrak m}_1 G(\tau) 
+ G (\overline{\mathfrak m}_1 \tau),
$$
where $\overline{\mathfrak m}_1=(-1)^{\dim L} \partial$.  
We can find such $G(\tau)$ by induction on dimension of $\tau$.  
In our case, we construct $G(\tau)$ using the gradient-like 
flow $\{\rho_t\}$.  
Set 
\[
{\mathfrak r}({\mathrm{Im}} \tau)=\bigcup_{t\leq 0} \rho_t(
{\mathrm{Im}} \tau).
\]
By the choice of our perturbations above, 
the closure of ${\mathfrak r}({\mathrm{Im}} \tau)$ can be 
triangulated in a compatible way with $\tau$ and $\Pi(\tau)$.  
Pick such a triangulation and then define $G(\tau)$ the corresponding 
singular chain.  
For the chain $G(\tau)$, we define $G(G(\tau))=0$.

Note that $\Pi:C_{(g)}(L) \to C_{(0)}(L)$ and 
$G:C_{(g)}(L) \to C_{(g)}(L)$ satisfy the conditions in Lemma 
3.2, hence $C_{(g)}(L)$ satisfies Condition 2.  
Therefore we can apply Theorem \ref{reduction} to reduce 
the filtered $A_{\infty}$-structure on $C^{\bullet}(L; 
\Lambda_{0,nov})$ to $CM^{\bullet}(f) \otimes \Lambda_{0,nov}$ 
and obtain a filtered $A_{\infty}$-algebra 
$(CM^{\bullet}(f) \otimes \Lambda_{0,nov},\{{\mathfrak m}'_k\})$, 
which is homotopy equivalent to $(C_{(g)}(L) \otimes \Lambda_{0,nov},
\{{\mathfrak m}_k\})$.  
Theorem \ref{Morse} is proved.   

In the proof of Theorem 3.3, we constructed the operator 
${\mathfrak m}'_k$ from ${\mathfrak m}_{\Gamma}$, 
$\Gamma \in G^+_{k+1}$.  
The geometric meaning of ${\mathfrak m}_{\Gamma}$ 
is as follows.  
Recall that $G(\tau)$ assigns the closure of the union of 
flow lines arriving at $\tau$.  
We assign $G$ to the interior edges.   
The interior vertices correspond to $J$-holomorphic discs, 
more precisely, bordered stable maps of genus 0.  
In order to describe the operation ${\mathfrak m}_{\Gamma}$,  
we need only rigid configuration of $\tau_i \in 
\chi_{\mathfrak T}(L)$ (the barycenters of $\tau_i$ are inputs), 
$J$-holomorphic discs, (broken) negative flow lines of $X$ and 
$W^s(q)$ ($q$ is the output).  
We choose the perturbation ${\mathfrak s}$ generically so that 
the moduli spaces of holomorphic discs and the flow $\{\rho_t\}$ 
are in general position so that 
the inner edges correspond to negative flow lines of $X$.  
Hence the ${\mathfrak m}_{\Gamma}$ is defined by using 
the configuration of pseudo-holomorphic discs and 
Morse negative gradient trajectories according to the decorated tree 
$\Gamma \in \cup_k G^+_{k+1}$.  

For a decorated tree $\Gamma \in G_{k+1}^+$, 
each edge is oriented in the direction 
from the $k$ input vertices to the root vertex.  
We denote by $v^{\pm}(e)$ the vertices such that 
$e$ is an oriented edge from $v^-(e)$ to $v^+(e)$.  
Consider 
the moduli space ${\mathcal M}_{\Gamma}(h;p_1, \dots. p_k, q)$ 
consisting of the configuration of the following 

\begin{itemize}
\item 
for each interior vertex $v \in \Gamma$, a bordered 
stable map $u_v$ representing the class $\beta_{\eta(v)}$ with 
$\ell (v)$ boundary marked points, 
where $\ell(v)$ is the valency of $v$, 
(we denote by $p(e,v)$ the marked point corresponding to the  
edge $e$ attached to $v$)
\item
the $i$-th input edge $e_i$ corresponds to a broken negative 
gradient flow line $\gamma_i$ starting from the critical point $p_i$ 
to $u_{v^+(e_i)}(p(e_i,v^+(e_i)))$, 
\item
the output edge corresponds to a broken negative gradient flow line 
$\gamma_0$ from $u_{v^-(e_0)}(p(e_0,v^-(e_0)))$ 
ending at the critical point $q$,
\item
an interior edge $e$ corresponds to a broken negative gradient flow line $\gamma_e$ from $u_{v^-(e)}(p(e,v^-(e)))$ to $u_{v^+(e)}(p(e,v^+(e)))$. 
\end{itemize}

Counting the weighted order of the moduli spaces of vitual 
dimension $0$, we get
$$
{\mathfrak m}_{\Gamma}(p_1 \otimes \dots \otimes p_k) 
= \sum_q  \# {\mathcal M}_{\Gamma}(h;p_1, \dots, p_k,q) \cdot 
e^{\sum_{v} \mu(\beta_{\eta(v)})/2}q
$$
and 
$${\mathfrak m}_k= \sum_{\Gamma \in G^+_{k+1}} 
T^{E(\Gamma)} {\mathfrak m}_{\Gamma}.
$$

\begin{figure}
\begin{center}
\include{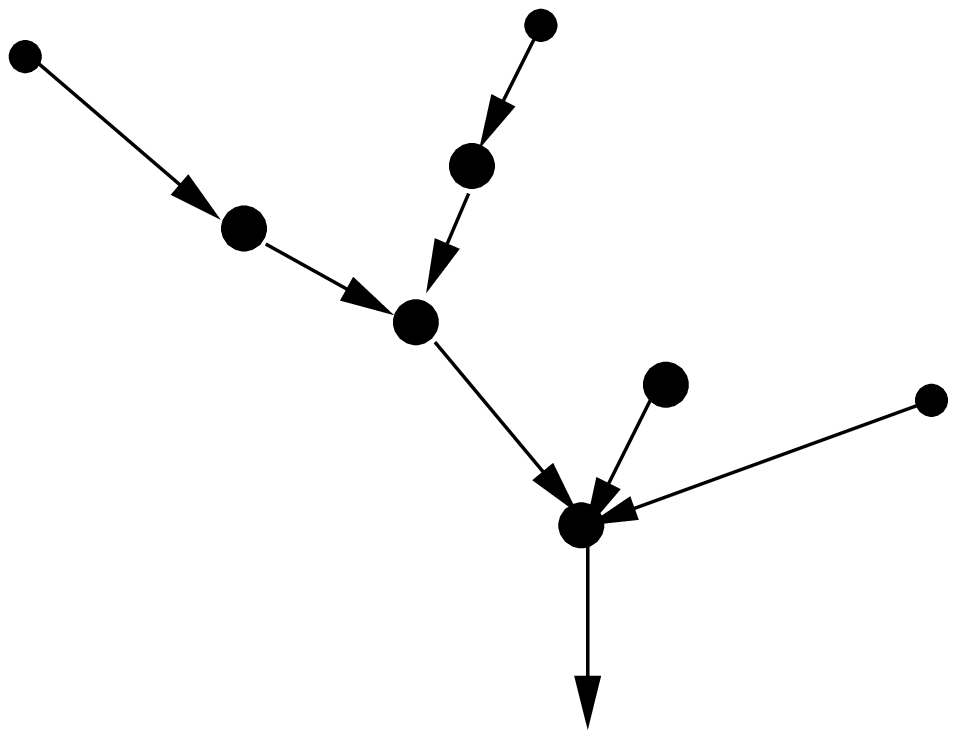}
\end{center}
\caption{}
\label{modelconfiguration}
\end{figure}

\begin{figure}
\begin{center}
\include{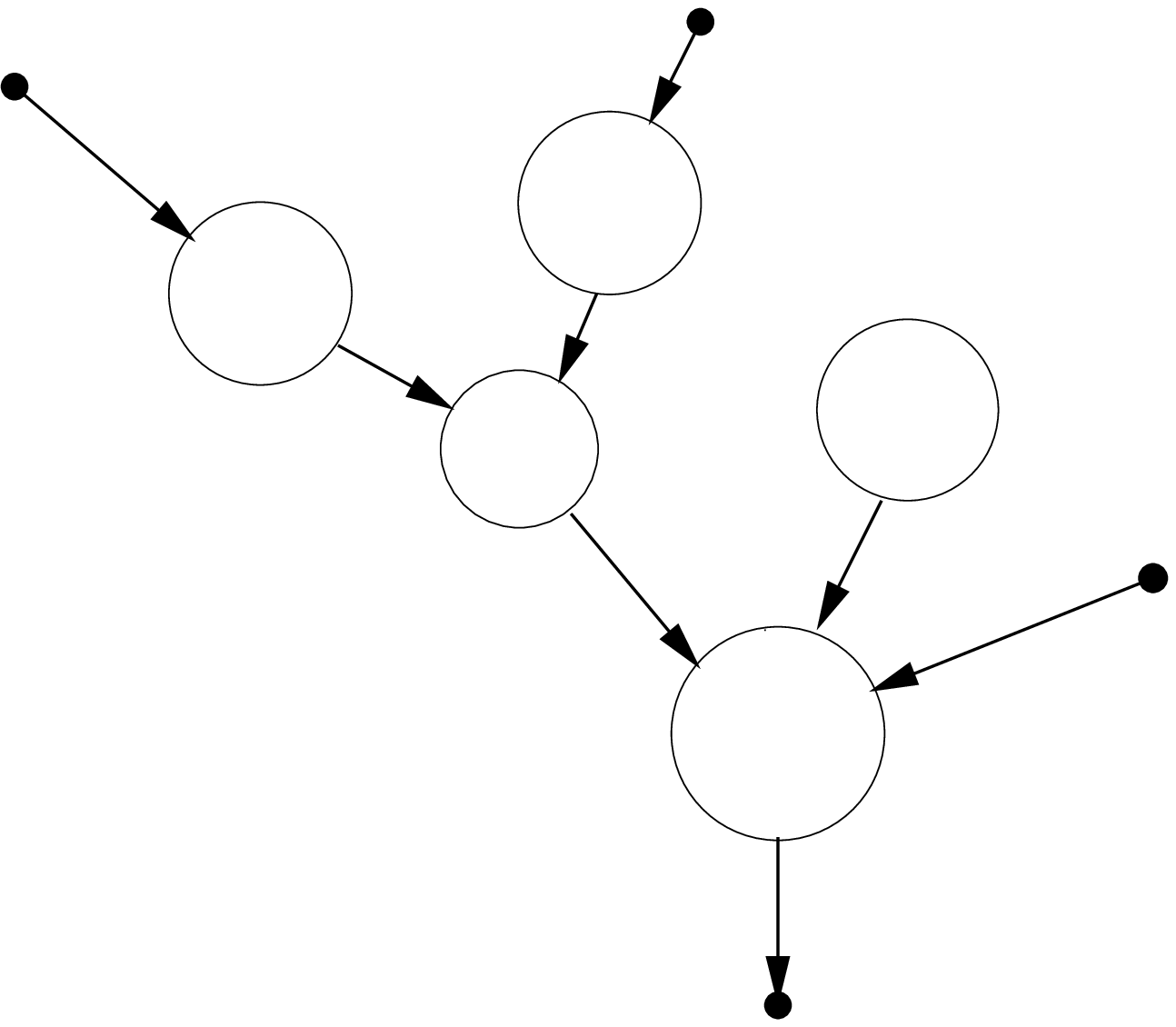}
\end{center}
\caption{}
\label{configuration}
\end{figure}

For example, 
we obtain the configuration as in Figure \ref{configuration} associated 
to the decorated planar tree $\Gamma$ 
with inputs $T(p),T(p'),T(p'') \in 
\chi_{\mathfrak T}(L)$  as in Figure \ref{modelconfiguration}.

This is essentially the configuration introduced in \cite{Fuk94}.  
Note that the first named author \cite{Fuk94} took multiple 
Morse functions to achieve transversality.  
Here we use one Morse function and apply the argument in section 3 
to squeeze the filtered 
$A_{\infty}$-algebra structure to the Morse complex.   
We emphasize that this becomes possible only after working out 
the chain level intersection theory in detail, which we explained 
in section 4.  
To find an appropriate perturbation of 
${\mathcal M}_{\Gamma}(h;p_1, \dots, p_k,q)$ directly without using 
the argument in section 4 (or section 30 in \cite{FOOO}) seems 
extremely difficult.  

The use of multiple Morse functions enables to construct the 
topological (or partial) filterd $A_{\infty}$-category of 
Morse functions on $L$ in the case that ${\mathfrak m}_0 = 0$.  
Note that, in a topological (or partial) filtered $A_{\infty}$-category 
${\mathcal A}$, 
the set $Ob_{\mathcal A}$ of objects is a topological space and 
the set $Mor_{\mathcal A}(a,b)$ of morphisms is defined 
for $(a,b)$ in an open dense subset of $Ob_{\mathcal A} \times 
Ob_{\mathcal A}$.  
When ${\mathcal A}$ is a filtered $A_{\infty}$-category, each 
object $a$ is equipped with the filtered $A_{\infty}$-algebra 
$Mor_{\mathcal A}(a,a)$.   
In our case, the filtered $A_{\infty}$-algebra on Morse complex 
$CM^{\bullet}(f)$ 
corresponds to the filtered $A_{\infty}$-algebra associated to 
the object $f$.  
Note that, in the construction of this paper and in Theorem 5.1, 
we do not need to assume that ${\mathfrak m}_0=0$ 
in our construction. 

For a relative spin pair $(L_0,L_1)$ of Lagrangian submanifolds, 
which intersect transversely, we obtain 
the filtered $A_{\infty}$-bimodule over the filtered 
$A_{\infty}$-algebras on $CM^{\bullet}(f_i) \otimes \Lambda_{0,nov}$, 
where $f_i:L_i \to \R$, 
$i=0,1$, are Morse functions.  
When $L_0$ and $L_1$ are of clean intersection, 
there exists a certain local system $\Theta$ on $L_0 \cap L_1$ 
and we can reduce the filtered $A_{\infty}$-bimodule structure on 
$C^{\bullet}(L_0 \cap L_1;\Theta) \otimes \Lambda_{0,nov}$ to 
$CM^{\bullet}(h;\Theta) 
\otimes \Lambda_{0,nov}$ 
over the filtered $A_{\infty}$-algebras on $CM^{\bullet}(f_i) \otimes 
\Lambda_{0,nov}$.  
Here $h$ is a Morse function on $L_0 \cap L_1$, which may be 
disconnected with various dimensions.  
For the canonical models of filtered $A_{\infty}$-bimodules, 
see \cite{FOOO}.  

\smallskip \noindent
{\bf Acknowledgement.}  We thank Otto van Koert for his kind 
instruction of making figures in this article.



\begin{thebibliography}{99}

\bibitem[1]{BC} P. Biran and O. Cornea, Quantum structures for 
Lagrangian submanifolds, preprint, arXiv:0708.4221.  

\bibitem[2]{Buk} L. Buhovsky, Multiplicative structures in 
Lagrangian Floer homology, preprint, math.SG/0608063.  

\bibitem[3]{Fl} A. Floer, Morse theory for Lagrangian intersections, 
J. Differential Geom. 28 (1988), 513-547.  

\bibitem[4]{Fuk94} K. Fukaya, Morse homotopy and its quantization, 
Geometric Topology (Athens, GA, 1993),  
AMS/IP Studies in Adv. Math. {\bf 2}, Part 1, 
Amer. Math. Soc., 1997, 409-440.    

\bibitem[5]{2000} K. Fukaya, Y.-G. Oh, H. Ohta and K. Ono, 
Lagrangian intersection Floer theory -anomaly and obstruction-, 
preprint 2000, available at http://www.math.kyoto-u.ac.jp/~fukaya

\bibitem[6]{FOOO} K. Fukaya, Y.-G. Oh, H. Ohta and K. Ono, 
Lagrangian intersection Floer theory -anomaly and obstruction-, 
revised and expanded version of \cite{2000}, preprint 2006.

\bibitem[7]{FOOO2} K. Fukaya, Y.-G. Oh, H. Ohta and K. Ono, 
Lagrangian Floer theory on compact toric manifolds, I, preprint 2008, 
arXiv:0802.1703, II, preprint 2008 arXiv:0810.5654.  

\bibitem[8]{FO} K. Fukaya and K. Ono, Arnold conjecture and 
Gromov-Witten invariant, Topology 38 (1999), 933-1048. 

\bibitem[9]{G-M} P. Griffiths and J. Morgan, Rational homotopy 
theory and differential forms, Progress in Math. {\bf 16}, 
Birkh\"auser, 1981.  

\bibitem[10]{Kad} T. V. Kadeishvili, The algebraic structure in 
the homology of an $A(\infty)$-algebra, 
Soobshch. Acad. Nauk Gruzin. SSR 108(1982), 249-252.  

\bibitem[11]{KS} M. Kontsevich and Y. Soibelman, 
Homological mirror symmetry and torus fibrations, 
Symplectic Geometry and Mirror Symmetry, edited by 
K. Fukaya, Y.-G. Oh, K. Ono, G. Tian, p.203-263, World Scientific 2001.  
\bibitem[12]{Oh} Y.-G. Oh, Floer cohomology of Lagrangian intersections 
and pseudo-holomorphic disks, I, Comm. Pure Appl. Math. 46 (1993), 
949-993.  

\bibitem[13]{Oh2} Y.-G. Oh, Relative Floer and quantum cohomology 
and the symplectic topology of Lagrangian submanifolds, 
Contact and symplectic geometry, 201-267, Publ. Newton Inst., {\bf 8}, 
Cambridge Univ. Press, 1996.  

\bibitem[14]{Smi} V. A. Smirnov, Simplicial and operad methods in 
algebraic topology, Translations of Mathematical Monographs 198, 
Amer. Math. Soc., 2001. 

\bibitem[15]{St} J. Stasheff, Homotopy associativity of H-spaces, 
I, II, Trans. Amer. Math. Soc. 108(1963), 275-312

\bibitem[16]{Su} D. Sullivan, Infinitesimal computations in topology, 
IHES Publ. Math. {\bf 47}(1977), 269-331.  
\end{thebibliography}
\end{document}